\documentclass[a4paper,11pt]{article}
\usepackage[a4paper, total={7in, 9.6in}]{geometry}
\usepackage[utf8]{inputenc} 
\usepackage[T1]{fontenc} 
\usepackage{amsmath}
\usepackage{amsthm}
\usepackage{verbatim}
\usepackage{mathtools} 
\usepackage{extarrows} 
\usepackage{amsthm}
\usepackage{xparse}
\usepackage{hyperref}
\usepackage{cleveref}
\usepackage{setspace}
\usepackage{mathpazo} 
\usepackage{amssymb}
\usepackage{tikz-cd}
\usepackage{float}
\usepackage{tabstackengine}
\setstackgap{L}{3.5mm}
\usepackage{graphicx}
\usepackage{tikz}
\usepackage{xcolor}
\usepackage{amssymb,latexsym} 
 
\numberwithin{equation}{section} 
 
\newtheorem{theorem}{Theorem}[section] 
\newtheorem{notation}[theorem]{Notation}
\newtheorem*{teo*}{Theorem}
\newtheorem{proposition}[theorem]{Proposition}

\newtheorem{lemma}[theorem]{Lemma} 
 
\theoremstyle{definition} 
 
\newtheorem{remark}[theorem]{Remark} 
\newtheorem{example}[theorem]{Example} 
\newtheorem{definition}[theorem]{Definition}

\usepackage[maxbibnames=99,style=alphabetic]{biblatex}
\usepackage{ytableau}
\usepackage{hyperref}
\usepackage[autostyle=true]{csquotes}

\begin{document}
\begin{center}
\textbf{\Large Cluster expansion formulas and perfect matchings for type B and C} 
\vspace{5mm}
 \\Azzurra Ciliberti\footnote{DIPARTIMENTO DI MATEMATICA “GUIDO CASTELNUOVO”, SAPIENZA UNIVERSIT\`{A} DI ROMA.\\
 \textit{Email address}: \textbf{azzurra.ciliberti@uniroma1.it}}

 \end{center}
 \vspace{0.5cm}
 \begin{abstract}
  \noindent 
Let $\mathbf{P}_{2n+2}$ be the regular polygon with $2n+2$ vertices, and let $\theta$ be the rotation of 180$^\circ$. Fomin and Zelevinsky proved that $\theta$-invariant triangulations of $\mathbf{P}_{2n+2}$ are in bijection with the clusters of cluster algebras of type $B_n$ and $C_n$. Furthermore, cluster variables correspond to the orbits of the action of $\theta$ on the diagonals of $\mathbf{P}_{2n+2}$. In this paper, we associate a labeled modified snake graph $\mathcal{G}_{ab}$ to each $\theta$-orbit $[a,b]$, and we get the cluster variables of type $B_n$ and $C_n$ which correspond to $[a,b]$ as perfect matching Laurent polynomials of $\mathcal{G}_{ab}$. This extends the work of Musiker for cluster algebras of type B and C to every seed.
\end{abstract}
\vspace{1cm}
\section*{Introduction}

Cluster algebras, discovered by Fomin and Zelevinsky and introduced in their seminal work \cite{CAI}, are commutative algebras with a special combinatorial structure. A \emph{cluster algebra} is a subalgebra of a field of rational functions in $n$ variables $u_1,\dots,u_n$ that is generated by the so called \emph{cluster variables}. Cluster variables are constructed recursively from an initial seed by a process called \emph{mutation}, and they are grouped into overlapping sets of constant cardinality $n$, the \emph{clusters}.

A remarkable result in the theory states that every cluster variable $x$ is a Laurent polynomial in the cluster variables $u_1,\dots,u_n$ of the initial cluster, i.e.,
\begin{equation}\label{cluster_exp}
    x=\frac{f(u_1,\dots,u_n)}{\displaystyle\prod_{i=1}^n u_i^{d_i}},
\end{equation}
where $f$ is a polynomial, and $d_1,\dots,d_n$ are non-negative integers. This is usually referred to as the \emph{Laurent phenomenon} \cite{CAI}, and the right hand side of equation \ref{cluster_exp} is the \emph{cluster expansion} of $x$ in the initial cluster variables.

Cluster algebras are related to a number of research areas including representation theory of finite dimensional algebras and Lie algebras, combinatorics, algebraic and hyperbolic geometry, dynamical systems, and knot theory.

Fomin, Shapiro and Thurston in \cite{FST,FT} initiate the study of cluster algebras coming from triangulations of surfaces with boundary and marked points. In their approach, cluster variables correspond to arcs in the surface, and clusters correspond to triangulations. Musiker and Schiffler in \cite{MS} give an expansion formula for the cluster variables in terms of perfect matchings of some labeled graphs that are recursively constructed from the surface by gluing together elementary pieces called tiles.

Let $\mathbf{P}_{2n+2}$ be the regular polygon with $2n+2$ vertices. Let $\theta$ be the rotation of $180^\circ$. Fomin and Zelevinsky show in \cite{CAII} that $\theta$-invariant triangulations of $\mathbf{P}_{2n+2}$ are in bijection with the clusters of cluster algebras of type $B_n$ and $C_n$. Furthermore, cluster variables correspond to the orbits of the action of $\theta$ on the diagonals of $\mathbf{P}_{2n+2}$. In \cite{ciliberti2024}, we consider cluster algebras of type $B$ and $C$ with principal coefficients in the initial seed. Using their model, we establish a formula that relates cluster variables of type $B_n$ and $C_n$ to cluster variables of type $A_n$.

In this note, we associate to each $\theta$-orbit $[a,b]$ of $\mathbf{P}_{2n+2}$ a labeled modified snake graph $\mathcal{G}_{ab}$ constructed by gluing together the snake graphs corresponding to particular diagonals obtained from those of $[a,b]$ by identifying some vertices of the polygon. Then we get the cluster expansion of the cluster variable $x_{ab}$ corresponding to $[a,b]$ in terms of perfect matchings of $\mathcal{G}_{ab}$. In particular, we present a combinatorial description of its $F$-polynomial and its $\bold{g}$-vector. This extends to every seed the work of Musiker \cite{M}, which provides cluster expansions for types $B$ and $C$ in terms of perfect matchings of modified snake graphs only for the initial bipartite seed. Another notable novelty is the connection with triangulations of polygons, which is not present in his approach. The main result of the paper is the following:

\begin{teo*}[\ref{theorem 3},\ref{theorem 4}]
 Let $T$ be a $\theta$-invariant triangulation of $\mathbf{P}_{2n+2}$. Let $\mathcal{A}=\mathcal{A}(T)^B$ (resp. $\mathcal{A}=\mathcal{A}(T)^C$) be the cluster algebra of type $B_n$ (resp. $C_n$) with principal coefficients in $T$. 
    Let $[a,b]$ be a $\theta$-orbit, and $x_{ab}$ the cluster variable of $\mathcal{A}$ corresponding to $[a,b]$. Let $F_{ab}$ and $\bold{g}_{ab}$ denote the $F$-polynomial and the $\bold{g}$-vector of $x_{ab}$, respectively. 
    Then $F_{ab}=F_{\mathcal{G}_{ab}}$ and $\bold{g}_{ab}=\bold{g}_{\mathcal{G}_{ab}}$, where $F_{\mathcal{G}_{ab}}$ is the perfect matching polynomial of $\mathcal{G}_{ab}$ and $\bold{g}_{\mathcal{G}_{ab}}$ is its $\bold{g}$-vector.  
\end{teo*}


Several other works in the literature use different techniques to study cluster expansion formulas for cluster algebras of type $B$ and $C$. In particular, Nakanishi and Stella provide in \cite{NS} a diagrammatic description of the $\bold{g}$-vectors of cluster algebras of type $B$ and $C$, while Reading studies them in \cite{R} using ring homomorphisms between cluster algebras of type $B$ and $C$, and cluster algebras of type $A$, induced by the fact that exchange matrices of type $B_n$ and $C_n$ ``dominate'' exchange matrices of type $A_n$. Moreover, a cluster algebra of type $B_n$ (resp. $C_n$) can be realized as a disk with one orbifold point of weight 2 (resp. $\frac{1}{2}$), and $n+1$ boundary marked points \cite{FeST}. In \cite{FeliksonTumarkin}, Felikson and Tumarkin compute $\bold{g}$-vectors for cluster algebras from orbifolds, including type $B$ and $C$, in terms of laminations on the orbifolds. In \cite{canakciTumarkin}, \text{\c{C}anak\c{c}\i} and Tumarkin introduce snake and band graphs associated to curves in a triangulated orbifold with orbifold points of weight $\frac{1}{2}$, including type $C$. Furthermore, a relation between skew-symmetric and skew-symmetrizable cluster algebras has been investigated in \cite{FST_finite_type,Dupont} via folding. Finally, in \cite{BanaianKelley}, Banaian and Kelley extend the snake graph construction of Musiker, Schiffler and Williams \cite{MSW11} to generalized cluster algebras from unpunctured orbifolds, including generalized cluster algebras of type $B$ and $C$.

The paper is organized as follows. In Section \ref{section_sg}, we recall from \cite{CSI} how to associate a snake graph to a diagonal in a triangulated polygon, and some basic notions about the combinatorics of these objects. In Section \ref{c_exp_formulas}, we report from \cite{ciliberti2024} cluster expansion formulas for cluster algebras of type $B$ and $C$ associated with $\theta$-invariant triangulations of the polygon.
Finally, in Section \ref{modified_sg}, we give the definition of labeled modified snake graph of a $\theta$-orbit, and prove the main result.

\section{Snake graphs from polygons}\label{section_sg}
Let $\mathbf{P}_{n+3}$ be a polygon with $n+3$ vertices and let $\Bar{T}=\{ \tau_1,\dots,\tau_n \}$ be a triangulation. Let $\gamma$ be a diagonal of $\mathbf{P}_{n+3}$ that is not in $\Bar{T}$. We choose an orientation on $\gamma$ such that $s$ is its starting point and $t$ its endpoint. Let

\begin{center}
    $s = p_0,p_1,p_2,\dots,p_{d+1} = t$
\end{center}
the intersections of $\gamma$ with $\Bar{T}$ in order of appearance, with $p_j \in \tau_{i_j}$. Let $\Delta_{j-1}$ and $\Delta_j$ be the two triangles of $\Bar{T}$ on each side of $\tau_{i_j}$. Let $G_j$ be the graph with 4 vertices and 5 edges, having the shape of a square with a diagonal, such that there is a bijection between the edges of $G_j$ and the 5 diagonals in the two triangles $\Delta_{j-1}$ and $\Delta_j$, where the diagonal in $G_j$ corresponds to the diagonal $\tau_{i_j}$. Moreover, this bijection must preserve the signed adjacency of the diagonals up to sign, that is, maintain the relative positioning of diagonals with respect to each other up to sign. 
 
 \begin{definition}
     The graph $G_j$ described above is called \emph{tile}.
 \end{definition}
 
 Given a planar embedding $\Tilde{G}_j$ of $G_j$ the relative orientation $\text{Rel}(\Tilde{G}_j,\Bar{T})$ of $\Tilde{G}_j$ with respect to $\Bar{T}$ is $+1$ (resp. $-1$) if its triangles agree (resp. disagree) in orientation with those of $\Bar{T}$.

 Diagonals $\tau_{i_j}$ and $\tau_{i_{j+1}}$ form two edges of the triangle $\Delta_j$ in $\Bar{T}$. The third edge of this triangle is labeled $\tau_{[j]}$.
Tiles $G_1, \dots , G_d$ in order from 1 to $d$ are glued together in the following way: $G_{j+1}$ is glued to $\Tilde{G}_j$, along the edge labeled $\tau_{[j]}$, choosing a planar embedding $\Tilde{G}_{j+1}$ for $G_{j+1}$ such that $\text{rel}(\Tilde{G}_{j+1},\Bar{T})\neq \text{rel}(\Tilde{G}_{j},\Bar{T})$, as in Figure \ref{gluing}. The resulting graph embedded in the plane is denoted by $\mathcal{G}_\gamma^\Delta$. The edges along which we glue two tiles are called \emph{internal}; the other ones are called \emph{external}.

\begin{figure}[H]
    \centering
   \begin{tikzpicture}[scale=2]
  
 
   \draw (-0.5,0) -- (-0.5,1) -- (0.5,1) -- (0.5,0) -- cycle;

\begin{scope}
    \draw (-0.5,0) --node[midway] {} (-0.5,1) --node[midway, above left,xshift=2mm] {} (0.5,1) -- node[midway,right,xshift=-1mm] {}(0.5,0) -- node[midway, below,yshift=1mm] {$\tau_{i_j}$} cycle;
    \draw (-0.5,1) -- node[midway,right,xshift=-1mm] {$\tau_{i_{j+1}}$}(0.5,0);

\end{scope}

\begin{scope}[xshift=-1cm]
    \draw (-0.5,0) -- (-0.5,1) -- node[midway, above,yshift=-1mm] {$\tau_{i_{j+i}}$}(0.5,1) -- node[midway, right,xshift=-1mm] {$\tau_{[j]}$}(0.5,0) --  cycle;
     \draw (-0.5,1) -- node[midway,right,xshift=-1mm] {$\tau_{i_{j}}$}(0.5,0);
\end{scope}

\end{tikzpicture}
\caption{Gluing tiles $\Tilde{G}_j$ and $\Tilde{G}_{j+1}$ along the edge labeled $\tau_{[j]}$.}
\label{gluing}
\end{figure}

\begin{definition}[{\cite[Definition 4.18]{CSI}}]\label{def_snake_graph}
The snake graph $\mathcal{G}_\gamma$ associated to $\gamma$ is obtained from $\mathcal{G}_\gamma^\Delta$ by removing the diagonal in each tile.    
\end{definition}
The edges of $\mathcal{G}_\gamma$ along which two tiles are glued are called \emph{internal}, while the remaining edges are referred to as \emph{external}.
\begin{definition}
A \emph{perfect matching} of a graph $\mathcal{G}$ is a subset $P$ of the edges of $\mathcal{G}$ such that each vertex of $\mathcal{G}$ is incident to exactly one edge of $P$.    
\end{definition}

\begin{definition}[{\cite[Definition 4.22]{CSI}}]\label{p-}
    Let $\gamma$ be a diagonal. The snake graph $\mathcal{G}_\gamma$ has precisely two perfect matchings which contain only boundary edges. If $\text{Rel}(\Tilde{G}_1,\Bar{T})=+1$ (resp. $-1$), $e_1$ and $e_2$ are defined to be the two edges of $\mathcal{G}_\gamma^\Delta$ which lie in the counterclockwise (resp. clockwise) direction from the diagonal of $\Tilde{G}_1$. Then $P_-=P_-(\mathcal{G}_\gamma^\Delta)$ is the unique matching which contains only boundary edges and does not contain edges $e_1$ or $e_2$. $P_-$ is called the \emph{minimal matching}. $P_+=P_+(\mathcal{G}_\gamma^\Delta)$, the \emph{maximal matching}, is the other matching with only boundary edges.
\end{definition}

Let $P_-\ominus P = (P_- \cup P) \setminus (P_- \cap P)$ be the symmetric difference of a perfect matching $P$ of the snake graph $\mathcal{G}_\gamma$.
By \cite[Theorem 5.1]{MS}, $P_-\ominus P$ is the set of boundary edges of a subgraph $\mathcal{G}_P$ of $\mathcal{G}_\gamma$, and $\mathcal{G}_P$ is a union of tiles
\begin{center}
    $\mathcal{G}_P=\displaystyle\bigcup_{i \in I} G_i$.
\end{center}
\begin{remark}
    The set $I$ depends on $P$.
\end{remark}
\begin{definition}[{\cite[Definition 4.24]{CSI}}]\label{def_h(P)_original}
Let $P$ be a perfect matching of $\mathcal{G}_\gamma$. The $height$ $monomial$ of $P$ is
\begin{center}
    $y(P):=\displaystyle\prod_{i \in I}y_i$.
\end{center}
Thus $y(P)$ is the product of all $y_i$ for which the tile $G_i$ lies inside $P_-\ominus P$.
\end{definition}

\begin{lemma}\label{lemma_ind_set}
    Let $\Tilde{I}=\{ i \mid \text{$(P_- \cup P)_{|G_i}$ contains an external edge of $\mathcal{G}_\gamma$ and $(P_- \cap P)_{|G_i}=\emptyset $}\}$. Then $\Tilde{I}=I$.
\end{lemma}

\begin{proof}
    If $i \in \Tilde{I}$, we have two cases to consider:
    \begin{itemize}
        \item [1)] $P_{|G_i}$ and $(P_-)_{|G_i}$ are both non-empty. In this case, since $P_-$ and $P$ are perfect matchings of $\mathcal{G}_\gamma$ which do not have any edges in common at the level of $G_i$, their union must include all external edges of $G_i$, so $i \in I$.
        \item [2)] Either $P_{|G_i}$ or $(P_-)_{|G_i}$, say $P_{|G_i}$, is empty. This means that two opposite edges of $G_i$ are internal in $\mathcal{G}_\gamma$. Since $P_-$ is a perfect matching of $\mathcal{G}_\gamma$ and $(P_-)_{|G_i}$ contains an external edge $e$ of $\mathcal{G}_\gamma$, $P_-$ must also contain the external edge opposite to $e$. Hence, $i \in I$.
    \end{itemize}
    Vice versa, if $i \in I$, then $(P_- \cup P)_{|G_i}$ contains all external edges of $G_i$ and $(P_- \cap P)_{|G_i} = \emptyset$.
\end{proof}

\begin{remark}\label{rmk_ind_set}
    It follows from Lemma \ref{lemma_ind_set}, that 
    $y(P)$ is the product of all $y_i$ such that $(P_- \cup P)_{|G_i}$ contains an external edge of $\mathcal{G}_\gamma$ and $(P_- \cap P)_{|G_i}=\emptyset $.
    
\end{remark}

\begin{definition}
    Let $\gamma$ be a diagonal which is not in $\Bar{T}$, and $\tau_{i_1}, \dots, \tau_{i_d}$ be the diagonals of $\Bar{T}$ crossed by $\gamma$. Then the \emph{perfect matching polynomial of} $\mathcal{G}_\gamma$ is
  \begin{center}
      $F_{\mathcal{G}_\gamma}:=\displaystyle\sum_{P}y(P)$,
  \end{center}
  where the sum is over all perfect matchings $P$ of $\mathcal{G}_\gamma$, and the $\bold{g}$-$vector$ is
  \begin{center}
      $\bold{g}_{\mathcal{G}_\gamma}:=\displaystyle\sum_{\tau_i \in P_-(\mathcal{G}_\gamma)}\bold{e}_i-\displaystyle\sum_{j=1}^d \bold{e}_{i_j}$,
  \end{center}
where $\{\bold{e}_1,\dots, \bold{e}_n\}$ is standard basis of $\mathbb{Z}^n$.
  The definition is extended to any diagonal by letting $F_{\mathcal{G}_\gamma}:=1$ and $\bold{g}_{\mathcal{G}_\gamma}:=\bold{e}_i$ if $\gamma=\tau_i \in \Bar{T}$, and $F_{\mathcal{G}_\gamma}:=1$ and $\bold{g}_{\mathcal{G}_\gamma}:=\bold{0}$ if $\gamma$ is a boundary edge of $\mathbf{P}_{n+3}$. 
\end{definition}

Let $\Bar{T}=\{\tau_1, \dots, \tau_n\}$ be a triangulation of $\mathbf{P}_{n+3}$. The \emph{adjacency matrix of $\Bar{T}$} is the skew-symmetric matrix $B(\Bar{T})=(b_{ij})$ such that $b_{ij}=1$ if $\tau_i$ and $\tau_j$ are two sides of a triangle of $\Bar{T}$ with $\tau_i$ following $\tau_j$ in counterclockwise order. The \emph{cluster algebra $\mathcal{A}(\Bar{T})$ of type $A_n$ with principal coefficients in $\Bar{T}$} is defined as the cluster algebra with principal coefficients in the initial seed whose exchange matrix is $B(\Bar{T})$ (see \cite[Example 6.6]{FST}). 

The following result is Theorem 3.1 of \cite{MS} restated in the case of polygons.
\begin{theorem}\label{f-poly g-vect of gamma}
  Let $\Bar{T}$ be a triangulation of $\mathbf{P}_{n+3}$, and let $\gamma$ be a diagonal which is not in $\Bar{T}$.  Then $F_{\mathcal{G}_\gamma}$ and $\bold{g}_{\mathcal{G}_\gamma}$ are the $F$-polynomial $F_\gamma$ and the $\bold{g}$-vector $\bold{g}_\gamma$ respectively of the cluster variable $x_\gamma$ of $\mathcal{A}(\Bar{T})$ which corresponds to $\gamma$. 
\end{theorem}

Given two distinct vertices $a$ and $b$ of $\mathbf{P}_{n+3}$, $(a,b)$ denotes the diagonal that connects them. We restate Definition 17.2 of \cite{FT} in the case of diagonals of a polygon as follows:

\begin{definition} \label{elem-laminate}
Let $\gamma=(a,b)$ be a diagonal of $\mathbf{P}_{n+3}$. The \emph{elementary lamination} associated to $\gamma$ is the segment $L_{\gamma}$ which 
 begins at $a'\in \mathbf{P}_{n+3}$ near $a$ in
the clockwise direction, and ends at $b' \in \mathbf{P}$ near $b$ in
the clockwise direction. 
If $\Bar{T}=\{ \tau_1,\dots,\tau_n \}$, then $L_i:=L_{\tau_i}$.
\end{definition}

 \begin{figure}[H]
        \centering
\begin{tikzpicture}[scale=0.7]
     \draw (90:3cm) -- (135:3cm) -- (180:3cm) -- (225:3cm) -- (270:3cm) -- (315:3cm) -- (360:3cm) -- (45:3cm) -- cycle;
                    \draw (90:3cm) -- node[near end, below, xshift=-1mm, yshift=-1mm] {} (180:3cm);
                    \draw (90:3cm) -- node[near end,below,xshift=1.2mm] {} (225:3cm);
                    \draw (90:3cm) -- node[near end, below,xshift=1mm] {} (270:3cm);
                    \draw (90:3cm) -- node[near end, above, xshift=-1mm] {} (315:3cm);
                    \draw (90:3cm) -- node[near end, above, xshift=-1mm] {} (360:3cm);

                    \draw[blue, line width=0.3mm] (80:3cm) -- node[near end, above, xshift=-1mm, yshift=-1mm] {$L_1$} (170:3cm);
                    \draw[blue, line width=0.3mm] (80:3cm) -- node[near end, above left,xshift=1.2mm] {$L_2$} (215:3cm);
                    \draw[blue, line width=0.3mm] (80:3cm) -- node[near end, above left,xshift=1mm] {$L_3$} (260:3cm);
                    \draw[blue, line width=0.3mm] (80:3cm) -- node[near end, below, xshift=-1.5mm] {$L_4$} (305:3cm);
                    \draw[blue, line width=0.3mm] (80:3cm) -- node[near end, below] {$L_5$} (350:3cm);

\end{tikzpicture}
        \caption{An octagon with the elementary lamination associated to each diagonal of the triangulation (in blue).}
        \label{lamination}
    \end{figure}

If $\gamma_1$ and $\gamma_2$ are two diagonals of $\mathbf{P}_{n+3}$ which cross each other, the exchange relation between cluster variables $x_{\gamma_1}$ and $x_{\gamma_2}$ of $\mathcal{A}(\Bar{T})$ can be phrased in terms of laminations of $\Bar{T}$ and diagonals of $\mathbf{P}_{n+3}$. The following Proposition is  Proposition 17.3 of \cite{FT} restated in the case of polygons. 

We denote by $(a,b)$ the diagonal of $\mathbf{P}_{n+3}$ which connects the vertices $a$ and $b$.

\begin{proposition} \label{up:skein1}
Let $\gamma_1=(a,b)$ and $\gamma_2=(c,d)$ be two diagonals of $\mathbf{P}_{n+3}$
which intersect. Then
\begin{equation} \label{u:skein-eq1}
x_{\gamma_1}x_{\gamma_2}=x_{(a,b)} x_{(c,d)} = \bold{y}^{\bold{d}_{ac,bd}} x_{(a,d)} ~x_{(b,c)}
 + \bold{y}^{\bold{d}_{ad,bc}}x_{(a,c)} ~x_{(b,d)},
\end{equation}
where
$\bold{d}_{ac,bd}$ (resp., $\bold{d}_{ad,bc}$) is the vector whose $i$-th coordinate is 1 if $L_i$ crosses both $(a,c)$ and $(b,d)$ (resp., $(a,d)$ and $(b,c)$); 0 otherwise.

    \begin{figure}[H]
        \centering
\begin{tikzpicture}[scale=0.7]
     \draw (90:3cm) -- (135:3cm) -- (180:3cm) -- (225:3cm) -- (270:3cm) -- (315:3cm) -- (360:3cm) -- (45:3cm) -- cycle;

                   \node at (135:3cm) [left] {$a$};
                   \node at (225:3cm) [left] {$c$};
                   \node at (315:3cm) [right] {$b$};
                   \node at (45:3cm) [right] {$d$};
                   \draw[red] (135:3cm) -- (315:3cm);
                   \draw[red] (225:3cm) -- (45:3cm);
                   \draw[green] (135:3cm) -- (225:3cm);
                   \draw[green] (315:3cm) -- (45:3cm);
                   \draw[yellow] (135:3cm) -- (45:3cm);
                   \draw[yellow] (225:3cm) -- (315:3cm);

\end{tikzpicture}
     \caption{The exchange relation between $x_{\gamma_1}$ and $x_{\gamma_2}$ viewed on the polygon.}   
    \end{figure}
    \end{proposition}

\section{Cluster expansion formulas for type B and C}\label{c_exp_formulas}
Let $\mathbf{P}_{2n+2}$ be the regular polygon with $2n+2$ vertices, and let $\theta$ be the rotation of $180^\circ$. Fomin and Zelevinsky show in \cite{CAII} that $\theta$-invariant triangulations of $\mathbf{P}_{2n+2}$ are in bijection with the clusters of a cluster algebra of type $B_n$ or $C_n$. Furthermore, cluster variables correspond to the orbits of the action of $\theta$ on the diagonals of $\mathbf{P}_{2n+2}$. 

 \begin{figure}[H]
                               \centering
                \begin{tikzpicture}[scale=0.5]
                    \draw (90:3cm) -- (135:3cm) -- (180:3cm) -- (225:3cm) -- (270:3cm) -- (315:3cm) -- (360:3cm) -- (45:3cm) -- cycle;
                    \draw (135:3cm) -- (225:3cm);;
                    \draw (45:3cm) -- (315:3cm);
                   \draw[] (270:3cm) -- node[midway, above left,xshift=1mm] {} (90:3cm);
                    \draw (135:3cm) -- node[midway, above left,xshift=1mm] {} (270:3cm);
                    \draw (315:3cm) -- node[midway, above left,xshift=1mm] {} (90:3cm);

                    \begin{scope}[xshift=8cm]
                      \draw (90:3cm) -- (135:3cm) -- (180:3cm) -- (225:3cm) -- (270:3cm) -- (315:3cm) -- (360:3cm) -- (45:3cm) -- cycle;
                    \draw (360:3cm) -- (270:3cm);
                    \draw (180:3cm) -- node[midway, above left,xshift=1mm] {} (90:3cm);
                    \draw (90:3cm) -- node[midway, above left,xshift=1mm] {} (270:3cm);
                     \draw (360:3cm) -- node[midway, above left,xshift=1mm] {} (90:3cm);
                    \draw (180:3cm) -- node[midway, above left,xshift=1mm] {} (270:3cm);

                    \end{scope}
                \end{tikzpicture}
                               \caption{Two $\theta$-invariant triangulations of $\mathbf{P}_8$.}
                               \label{ex_f_poly}
                           \end{figure}

Each $\theta$-invariant triangulation $T$ of $\mathbf{P}_{2n+2}$ has exactly one diameter $d$. After choosing an orientation of $d$, in \cite{ciliberti2024} we define cluster algebras $\mathcal{A}^B(T)$ of finite type $B_n$ and $\mathcal{A}^C(T)$ of finite type $C_n$, with principal coefficients in $T$. We also define a simple operation on the diagonals of $\mathbf{P}_{2n+2}$ which allows us to relate cluster variables of type $B_n$ and $C_n$ with those of type $A_n$.

\begin{definition}[{\cite[Definition 3.1]{ciliberti2024}}]\label{def_restriction}
    Let $\mathcal{D}$ be a set of diagonals of $\mathbf{P}_{2n+2}$. The \emph{restriction of $\mathcal{D}$}, denoted by $\text{Res}(\mathcal{D})$, is the set of diagonals of $\mathbf{P}_{n+3}$ obtained from those of $\mathcal{D}$ identifying all the vertices which lie on the right of $d$. 
\end{definition}

\begin{remark}
    The restriction depends on $d$, but not on the other diagonals of $T$.
\end{remark}

Let $D=\mathrm{diag}(1,\dots,1,2)\in \mathbb{Z}^n \times \mathbb{Z}^n$ be the diagonal matrix with diagonal entries $1,\dots,1,2$. Then $\mathcal{A}^B(T)$ (resp. $\mathcal{A}^C(T)$) is defined as the cluster algebra with principal coefficients in the initial seed whose exchange matrix is $DB(\Bar{T})$ (resp. $B(\Bar{T})D$). In \cite{ciliberti2024}, we provide a formula that expresses each cluster variable of type $B_n$ and $C_n$, in $\mathcal{A}^B(T)$ and $\mathcal{A}^C(T)$ respectively, in terms of cluster variables of type $A_n$ in $\mathcal{A}(\Bar{T})$, where $\Bar{T}=\text{Res}(T)$. In particular, denoting by $x_{ab}$ the cluster variable which corresponds to the $\theta$-orbit $[a,b]$ of the diagonal $(a,b)$, we combinatorially describe the $F$-polynomial $F_{ab}$ and the $\bold{g}$-vector $\bold{g}_{ab}$ of $x_{ab}$. In the rest of the section, we recall these results from \cite{ciliberti2024}.

In the following $T=\{ \tau_1, \dots, \tau_n=d, \dots, \tau_{2n-1} \}$ is a $\theta$-invariant triangulation of $\mathbf{P}_{2n+2}$ with oriented diameter $d$, and we assume that $\theta(\tau_i)=\tau_{2n-i}$. Moreover, we adopt the following notation:
\begin{notation}
\begin{itemize}
\item[i)] We use the label $\ast$ for the vertex of $\mathbf{P}_{n+3}$ obtained by identifying the vertices of $\mathbf{P}_{2n+2}$ which lie on the right of $d$.
    \item[ii)] For a vertex $a$ of $\mathbf{P}_{2n+2}$,  $\Bar{a}:=\theta(a)$.
    \end{itemize}
    
\end{notation}
\subsection{Type B}

\begin{definition}[{\cite[Definition 3.3]{ciliberti2024}}]\label{def_type_B}
Let $[a,b] \not \subset T$ be a $\theta$-orbit of $\mathbf{P}_{2n+2}$. Let $D=\mathrm{diag}(1,\dots,1,2)\in \mathbb{Z}^n \times \mathbb{Z}^n$ be the diagonal matrix with diagonal entries $1,\dots,1,2$.
\begin{itemize}
    \item[i)] If $\text{Res}([a,b])=\{ \gamma \}$ (as in Figure \ref{type B 1}), 
  \begin{equation}
      F_{ab}^B:=F_\gamma,
       \end{equation}
\begin{equation}
 \hspace{3cm}\bold{g}_{ab}^B:=\begin{cases}
          \text{$D\bold{g}_\gamma$ if $\gamma$ does not cross $\tau_n=d$;}\\
          \text{$D\bold{g}_\gamma+\bold{e}_n$ if $\gamma$ crosses $\tau_n=d$}.
      \end{cases}
      \end{equation}
  
  \item[ii)] Otherwise, $\text{Res}([a,b])=\{\gamma_1,\gamma_2\}$ (as in Figure \ref{type B}), and
\begin{equation}\label{def_F^1}
    F_{ab}^B:=F_{\gamma_1}F_{\gamma_2}- \bold{y}^{\bold{d}_{\gamma_1,\gamma_2}}F_{(a,\Bar{b})},
    \end{equation}
\begin{equation}\label{def_g^1}
    \hspace{-0.1cm}\bold{g}_{ab}^B:=D(\bold{g}_{\gamma_1}+\bold{g}_{\gamma_2}+\bold{e}_n).
\end{equation}

The definition is extended to any $\theta$-orbit by letting $F_{ab}^B:=1$ and $\bold{g}_{ab}^B:=\bold{e}_i$ if $[a,b]=\{ \tau_i,\tau_{2n-i} \} \in T$, and $F_{ab}^B:=1$ and $\bold{g}_{ab}^B:=\bold{0}$ if $(a,b)$ is a boundary edge of $\mathbf{P}_{2n+2}$.
\end{itemize}
   \begin{figure}[H]
                               \centering
                \begin{tikzpicture}[scale=0.5]
                    \draw (90:3cm) -- (120:3cm) -- (150:3cm) -- (180:3cm) -- (210:3cm) -- (240:3cm) -- (270:3cm) -- (300:3cm) -- (330:3cm) -- (360:3cm) -- (30:3cm) -- (60:3cm) --  cycle;
                  
                    \draw[-{Latex[length=2mm]}] (270:3cm) -- node[midway, above left,xshift=1mm] {} (90:3cm);
                  
                     \draw[red, line width=0.3mm] (150:3cm) -- (330:3cm);

                    \node at (330:3cm) [right] {$\Bar{a}$};

\node at (150:3cm) [left] {$a$};

                    \begin{scope}[xshift=8cm]
                      \draw (90:3cm) -- (120:3cm) -- (150:3cm) -- (180:3cm) -- (210:3cm) -- (240:3cm) -- (270:3cm) -- (360:3cm)  --  cycle;
                    \draw (90:3cm) -- node[midway, above left,xshift=1mm] {} (270:3cm);
                    \draw[red, line width=0.3mm] (150:3cm) -- node[midway, above right, xshift=-2mm, yshift=-1mm] {$\textcolor{black}{\gamma}$} (360:3cm); 
                  
                    \node at (360:3cm) [right] {$\ast$};

\node at (150:3cm) [left] {$a$};
                      
                    \end{scope}
                    \begin{scope}[yshift=-9cm]
                         \draw (90:3cm) -- (120:3cm) -- (150:3cm) -- (180:3cm) -- (210:3cm) -- (240:3cm) -- (270:3cm) -- (300:3cm) -- (330:3cm) -- (360:3cm) -- (30:3cm) -- (60:3cm) --  cycle;
                  
                    \draw[-{Latex[length=2mm]}] (270:3cm) -- node[midway, above left,xshift=1mm] {} (90:3cm);
                  
                     \draw[red, line width=0.3mm] (150:3cm) -- (270:3cm); 
                   
 \draw[red, line width=0.3mm] (90:3cm) -- (330:3cm);

                    \node at (90:3cm) [above] {$\Bar{b}$};
                    \node at (330:3cm) [right] {$\Bar{a}$};

\node at (150:3cm) [left] {$a$};
                    \node at (270:3cm) [below] {$b$};

                    \begin{scope}[xshift=8cm]
                      \draw (90:3cm) -- (120:3cm) -- (150:3cm) -- (180:3cm) -- (210:3cm) -- (240:3cm) -- (270:3cm) -- (360:3cm)  --  cycle;
                    \draw (90:3cm) -- node[midway, above left,xshift=1mm] {} (270:3cm);
                    \draw[red, line width=0.3mm] (150:3cm) -- node[midway, above right, xshift=-2mm, yshift=-1mm] {$\textcolor{black}{\gamma}$} (270:3cm); 
                  
                    \node at (360:3cm) [right] {$\ast$};
 \node at (270:3cm) [below] {$b$};
\node at (150:3cm) [left] {$a$};
                      
                    \end{scope}
                    \end{scope}
                \end{tikzpicture}
                               \caption{On the left, two $\theta$-orbits $[a,\Bar{a}]$ and $[a,b]$. On the right, their restrictions.}
                               \label{type B 1}
                           \end{figure}  

  \begin{figure}[H]
                               \centering
                \begin{tikzpicture}[scale=0.5]
                    \draw (90:3cm) -- (120:3cm) -- (150:3cm) -- (180:3cm) -- (210:3cm) -- (240:3cm) -- (270:3cm) -- (300:3cm) -- (330:3cm) -- (360:3cm) -- (30:3cm) -- (60:3cm) --  cycle;
                    \draw[-{Latex[length=2mm]}] (270:3cm) -- node[midway, above left,xshift=1mm] {} (90:3cm);
                     
                     \draw[red, line width=0.3mm] (150:3cm) -- (30:3cm); 
                    \draw[red, line width=0.3mm] (210:3cm) -- (330:3cm);

                    \node at (210:3cm) [left] {$\Bar{b}$};
                    \node at (330:3cm) [right] {$\Bar{a}$};

\node at (150:3cm) [left] {$a$};
                    \node at (30:3cm) [right] {$b$};

                    \begin{scope}[xshift=8cm]
                      \draw (90:3cm) -- (120:3cm) -- (150:3cm) -- (180:3cm) -- (210:3cm) -- (240:3cm) -- (270:3cm) -- (360:3cm)  --  cycle;
                    \draw (90:3cm) -- node[midway, above left,xshift=1mm] {} (270:3cm);
                    \draw[red, line width=0.3mm] (150:3cm) -- node[midway, above right, xshift=-2mm, yshift=-1mm] {$\textcolor{black}{\gamma_1}$} (360:3cm); 
                     \draw[cyan, line width=0.3mm] (150:3cm) -- (210:3cm);
\draw[red, line width=0.3mm] (210:3cm) -- node[midway, above right, xshift=-1mm] {$\textcolor{black}{\gamma_2}$} (360:3cm);
                  
                    \node at (210:3cm) [left] {$\Bar{b}$};
                    \node at (360:3cm) [right] {$\ast$};

\node at (150:3cm) [left] {$a$};
                      
                    \end{scope}
                \end{tikzpicture}
                               \caption{On the left, a $\theta$-orbit $[a,b]$. On the right, its restriction in red and the diagonal $(a,\Bar{b})$ in blue.}
                               \label{type B}
                           \end{figure}

\end{definition}

\begin{theorem}[{\cite[Theorem 3.4]{ciliberti2024}}]\label{theorem1}
    Let $\mathcal{A}=\mathcal{A}(T)^B$ be the cluster algebra of type $B_n$ with principal coefficients in $T$. 
    Let $[a,b]$ be a $\theta$-orbit, and $x_{ab}$ the cluster variable of $\mathcal{A}$ which corresponds to $[a,b]$. Let $F_{ab}$ and $\bold{g}_{ab}$ be the $F$-polynomial and the $\bold{g}$-vector of $x_{ab}$, respectively. 
    Then $F_{ab}=F_{ab}^B$ and $\bold{g}_{ab}=\bold{g}_{ab}^B$.
\end{theorem}

\subsection{Type C}
\begin{definition}[{\cite[Definition 3.9]{ciliberti2024}}]
    Let $T$ be a $\theta$-invariant triangulation of $\mathbf{P}_{2n+2}$, and let $\mathcal{O} \not \subset T$ be a $\theta$-orbit. The \emph{rotated restriction of $[a,b]$}, denoted by $\tilde{\text{Res}}([a,b])$, is defined as follows. 
    
    \begin{itemize}
        \item[i)] If $\mathcal{O}=[a,\bar{a}]$ is a diameter, so $\text{Res}([a,\bar{a}])=\{\gamma\}$, then $\tilde{\text{Res}}([a,\bar{a}]):=\{\tilde\gamma_1, \tilde\gamma_2\}$, where $\tilde\gamma_1=\gamma$ and $\tilde\gamma_2$, if it exists, is the diagonal of $\mathbf{P}_{n+3}$ which intersects the same diagonals of $T$ as $\gamma$ but $d$. If there is no such diagonal, $\tilde{\text{Res}}([a,\bar{a}]):=\{\tilde\gamma_1\}$.

        \item[ii)] If $\mathcal{O}=[a,b]$ is a pair of diagonals which do not cross $d$, then $\tilde{\text{Res}}([a,b]):=\text{Res}([a,b])$.

        \item[iii)] If $\mathcal{O}=[a,b]$ is a pair of diagonals which cross $d$, then $\text{Res}([a,b])=\{\gamma_1,\gamma_2\}$, where $\gamma_1$ and $\gamma_2$ are two diagonals of $\mathbf{P}_{n+3}$ that share the right endpoint, and such that $\gamma_2$ is obtained from $\gamma_1$ by rotating counterclockwise (resp. clockwise) its left endpoint if $\tau_{n-1}$ is counterclockwise (resp. clockwise) from $\tau_n$. Then $\tilde{\text{Res}}([a,b]):=\{\tilde\gamma_1, \tilde\gamma_2\}$, where  $\tilde\gamma_1=\gamma_1$ and $\tilde\gamma_2$, if it exists, is the diagonal of $\mathbf{P}_{n+3}$ which intersects the same diagonals of $T$ as $\gamma_2$ but the diameter. If there is no such diagonal, $\tilde{\text{Res}}([a,b]):=\{\tilde\gamma_1\}$.
      
    \end{itemize}
 
\end{definition}

\begin{definition}[{\cite[Definition 3.10]{ciliberti2024}}]
    Let $v \in \mathbb{Z}_{\geq 0}^{2n-1}$. The $rotated$ $restriction$ $of$ $v$, denoted by $\tilde{\text{Res}}(v)$, is the vector of the first $n$ coordinates of $v$, with the $n$-th one divided by 2.
\end{definition}

\begin{definition}[{\cite[Definition 3.11]{ciliberti2024}}]\label{def_type_C}
  Let $\mathcal{O} \not \subset T$ be a $\theta$-orbit of $\mathbf{P}_{2n+2}$.

  \begin{itemize}
      \item If $\mathcal{O}=[a,b]$ and $\tilde{\text{Res}}([a,b])=\{\tilde\gamma\}$,
      \begin{equation}
      F_{ab}^C:=F_{\Tilde{\gamma}},
\end{equation}
\begin{equation}
      \bold{g}_{ab}^C:=\begin{cases}
          \text{$\bold{g}_{\Tilde{\gamma}}+\bold{e}_i$ \hspace{0.5cm} if $\tau_i$ and $\tau_n=d$ are two different sides of a triangle of $T$},\\ \hspace{1.8cm}\text{$\tau_i$ is clockwise from $\tau_n$, and $\tilde\gamma$ crosses $\tau_n=d$;}\\
          \text{$\bold{g}_{\tilde\gamma}$ \hspace{1.3cm}otherwise.}
      \end{cases}
  \end{equation}
\item If $\mathcal{O}=[a,\Bar{a}]$ is a diameter, $\tilde{\text{Res}}([a,\Bar{a}])=\{\tilde\gamma_1, \tilde\gamma_2\}$, and it follows from the definition that there are uniquely determined two $\theta$-orbits $[a,\Bar{c}]$ and $[a,\Bar{b}]$, such that $\tilde{\text{Res}}([a,\Bar{c}])=\{\tilde\gamma_1\}$ and $\tilde{\text{Res}}([a,\Bar{b}])=\{\tilde\gamma_2\}$. A possible situation is represented in Figure \ref{type$C$-1}.

    Then
      \begin{equation}
        F_{a\Bar{a}}^C:=F_{\tilde\gamma_1}F_{\tilde\gamma_2}- \bold{y}^{\Tilde{\text{Res}}(\bold{d}_{a\ast,c\Bar{b}}+\bold{d}_{a\Bar{b},b\ast})}F_{(a,b)}F_{(a,c)},
\end{equation}
\begin{equation}\label{eq_g_vect_1}
      \bold{g}_{a\Bar{a}}^C:=\begin{cases}
          \text{$\bold{g}_{\tilde\gamma_1}+\bold{g}_{\tilde\gamma_2}+\bold{e}_i-\bold{g}_{(\Bar{b},c)}$  \hspace{0.6cm}if $\tau_i$ and $\tau_n$ are two different sides of a triangle of $T$},\\ \hspace{4.2cm}\text{and $\tau_i$ is clockwise from $\tau_n$;}\\
          \text{$\bold{g}_{\tilde\gamma_1}+\bold{g}_{\tilde\gamma_2}$ \hspace{2.6cm}otherwise.}
      \end{cases}
  \end{equation}
 \item If $\mathcal{O}=[a,b]$ is a pair of diagonals which cross $d$, and $\tilde{\text{Res}}([a,b])=\{\tilde\gamma_1,\tilde\gamma_2\}$, it follows from the definition that there are uniquely determined two $\theta$-orbits $[a,d]$ and $[b,c]$, such that $\tilde{\text{Res}}([a,d])=\{\tilde\gamma_1\}$ and $\tilde{\text{Res}}([b,c])=\{\tilde\gamma_2\}$. A possible situation is represented in Figure \ref{type$C$-2}.

 Then
      \begin{equation}
        F_{ab}^C:=F_{\tilde\gamma_1}F_{\tilde\gamma_2}- \bold{y}^{\Tilde{\text{Res}}(\bold{d}_{\Bar{b}\ast,\Bar{d}\Bar{c}}+\bold{d}_{a\Bar{c},c\ast})}F_{(a,c)}F_{(\Bar{b},\Bar{c})},
\end{equation}
\begin{equation}\label{eq_g_vect_2}
      \bold{g}_{ab}^C:=\begin{cases}
          \text{$\bold{g}_{\tilde\gamma_1}+\bold{g}_{\tilde\gamma_2}+\bold{e}_i-\bold{g}_{(\Bar{c},\Bar{d})}$  \hspace{0.6cm}if $\tau_i$ and $\tau_n$ are two different sides of a triangle of $T$},\\ \hspace{4.2cm}\text{and $\tau_i$ is clockwise from $\tau_n$;}\\
          \text{$\bold{g}_{\tilde\gamma_1}+\bold{g}_{\tilde\gamma_2}$ \hspace{2.6cm}otherwise.}
      \end{cases}
  \end{equation}    
  \end{itemize}
The definition is extended to any $\theta$-orbit by letting $F_{ab}^C:=1$ and $\bold{g}_{ab}^C:=\bold{e}_i$ if $[a,b]=\{ \tau_i,\tau_{2n-i} \} \in T$, and $F_{ab}^C:=1$ and $\bold{g}_{ab}^C:=\bold{0}$ if $(a,b)$ is a boundary edge of $\mathbf{P}_{2n+2}$.  
\end{definition}
      \begin{figure}[H]
                               \centering
                \begin{tikzpicture}[scale=0.5]
                    \draw (90:3cm) -- (120:3cm) -- (150:3cm) -- (180:3cm) -- (210:3cm) -- (240:3cm) -- (270:3cm) -- (300:3cm) -- (330:3cm) -- (360:3cm) -- (30:3cm) -- (60:3cm) --  cycle;
                    \draw[-{Latex[length=2mm]}] (270:3cm) -- node[midway, above left,xshift=1mm] {} (90:3cm);
                     
                    \draw[cyan, line width=0.3mm] (150:3cm) -- (270:3cm);
                     \draw[red, line width=0.3mm] (150:3cm) -- (60:3cm); 
                     \draw[red, line width=0.3mm] (330:3cm) -- (240:3cm);
                     \draw[cyan, line width=0.3mm] (330:3cm) -- (90:3cm);
                    \draw[black, line width=0.3mm] (150:3cm) -- (330:3cm);

                    \node at (90:3cm) [above,yshift=-1] {$b$};
                    \node at (270:3cm) [below,yshift=1] {$\Bar{b}$};

\node at (150:3cm) [left] {$a$};
                    \node at (60:3cm) [right] {$\Bar{c}$};
                    \node at (240:3cm) [left] {$c$};
                    \node at (330:3cm) [right] {$\Bar{a}$};

                    \begin{scope}[xshift=8cm]
                      \draw (90:3cm) -- (120:3cm) -- (150:3cm) -- (180:3cm) -- (210:3cm) -- (240:3cm) -- (270:3cm) -- (360:3cm)  --  cycle;
                    \draw (90:3cm) -- node[midway, above left,xshift=1mm] {} (270:3cm);
                     \draw[green, line width=0.3mm] (150:3cm) -- (90:3cm);
                   \draw[yellow, line width=0.3mm] (150:3cm) -- (240:3cm);
\draw[cyan, line width=0.3mm] (150:3cm) -- node[midway, above right, xshift=-2mm, yshift=-1mm] {$\textcolor{black}{\tilde\gamma_2}$}(270:3cm);
\draw[red, line width=0.3mm] (150:3cm) -- node[midway, above right, xshift=-2mm, yshift=-1mm] {$\textcolor{black}{\tilde\gamma_1}$}(360:3cm);

                    \node at (360:3cm) [right] {$\ast$};
 \node at (240:3cm) [left] {$c$};
\node at (150:3cm) [left] {$a$};
\node at (90:3cm) [above,yshift=-1] {$b$};
                    \node at (270:3cm) [below,yshift=1] {$\Bar{b}$};

                    \end{scope}
                \end{tikzpicture}
                               \caption{On the left, the $\theta$-orbits $[a,\Bar{a}]$, $[a,\Bar{c}]$, $[a,\Bar{b}]$. On the right, their rotated restrictions, and the diagonals $(a,b)$ and $(a,c)$.}
                               \label{type$C$-1}
                           \end{figure}  
     \begin{figure}[H]
                               \centering
                \begin{tikzpicture}[scale=0.5]
                    \draw (90:3cm) -- (120:3cm) -- (150:3cm) -- (180:3cm) -- (210:3cm) -- (240:3cm) -- (270:3cm) -- (300:3cm) -- (330:3cm) -- (360:3cm) -- (30:3cm) -- (60:3cm) --  cycle;
                    \draw[-{Latex[length=2mm]}] (270:3cm) -- node[midway, above left,xshift=1mm] {} (90:3cm);
                     
                    \draw[cyan, line width=0.3mm] (180:3cm) -- (270:3cm);
                     \draw[red, line width=0.3mm] (150:3cm) -- (60:3cm); 
                     \draw[red, line width=0.3mm] (330:3cm) -- (240:3cm);
                     \draw[cyan, line width=0.3mm] (360:3cm) -- (90:3cm);
                    \draw[black, line width=0.3mm] (150:3cm) -- (360:3cm);
                    \draw[black, line width=0.3mm] (180:3cm) -- (330:3cm);

                    \node at (90:3cm) [above,yshift=-1] {$c$};
                    \node at (270:3cm) [below,yshift=1] {$\Bar{c}$};
                    \node at (180:3cm) [left] {$\Bar{b}$};
                    \node at (360:3cm) [right] {$b$};
                    \node at (330:3cm) [right] {$\Bar{a}$};

\node at (150:3cm) [left] {$a$};
                    \node at (300:3cm) [right,xshift=-1] {$\hat{d}$};
                    \node at (60:3cm) [right] {$d$};
                    \node at (240:3cm) [left] {$\Bar{d}$};
                    \node at (30:3cm) [right] {$\hat{a}$};

                    \begin{scope}[xshift=8cm]
                      \draw (90:3cm) -- (120:3cm) -- (150:3cm) -- (180:3cm) -- (210:3cm) -- (240:3cm) -- (270:3cm) -- (360:3cm)  --  cycle;
                    \draw (90:3cm) -- node[midway, above left,xshift=1mm] {} (270:3cm);
                     \draw[green, line width=0.3mm] (150:3cm) -- (90:3cm);
                    \draw[yellow, line width=0.3mm] (180:3cm) -- (240:3cm);
\draw[cyan, line width=0.3mm] (180:3cm) -- node[midway, above right, xshift=-2mm, yshift=-1mm] {$\textcolor{black}{\tilde\gamma_2}$}(270:3cm);
\draw[red, line width=0.3mm] (150:3cm) -- node[midway, above right, xshift=-2mm, yshift=-1mm] {$\textcolor{black}{\tilde\gamma_1}$}(360:3cm);

                    \node at (360:3cm) [right] {$\ast$};
 \node at (240:3cm) [left] {$\Bar{d}$};
\node at (150:3cm) [left] {$a$};
\node at (90:3cm) [above,yshift=-1] {$c$};
                    \node at (270:3cm) [below,yshift=1] {$\Bar{c}$};
                     \node at (180:3cm) [left] {$\Bar{b}$};
                      
                    \end{scope}
                \end{tikzpicture}
                               \caption{On the left, the $\theta$-orbits $[a,b]$, $[a,d]$, $[b,c]$. On the right, their rotated restrictions, and the diagonals $(a,c)$ and $(\Bar{b},\Bar{d})$.}
                               \label{type$C$-2}
                           \end{figure} 
\begin{remark}\label{rmk_coeff_typeC}
We have that $\Tilde{\text{Res}}(\bold{d}_{a\ast,c\Bar{b}}+\bold{d}_{a\Bar{b},b\ast})$ is either equal to $\bold{d}_{a\ast,c\Bar{b}}$ or equal to $\bold{d}_{a\Bar{b},b\ast}$. Indeed, if $L_i$ and $L_j$, $i,j \neq n$, are the elementary lamination of two diagonals $\tau_i$ and $\tau_j$ of $T$ such that $L_i$ crosses both $(a,\ast)$ and $(c,\Bar{b})$, and $L_j$ crosses both $(a,\Bar{b})$ and $(b,\ast)$, then $L_i$ crosses $L_j$, so $\tau_i$ crosses $\tau_j$. Furthermore, if $L_n$ crosses both $(a,\ast)$ and $(c,\Bar{b})$, it also crosses $(a,\Bar{b})$ and $(b,\ast)$. Similarly, $\Tilde{\text{Res}}(\bold{d}_{\Bar{b}\ast,\Bar{d}\Bar{c}}+\bold{d}_{a\Bar{c},c\ast})$ is either equal to $\bold{d}_{\Bar{b}\ast,\Bar{d}\Bar{c}}$ or equal to $\bold{d}_{a\Bar{c},c\ast}$.   
\end{remark}

\begin{remark}
    We observe that $(\Bar{b},c)$ in \ref{eq_g_vect_1} and $(\Bar{c},\Bar{d})$ in \ref{eq_g_vect_2} are either diagonals of $\Bar{T}$ or boundary edges, since $\tilde{\text{Res}}([a,\Bar{c}])=\{\tilde\gamma_1\}$ and $\tilde{\text{Res}}([a,d])=\{\tilde\gamma_1\}$ respectively. Remember that, if $(a,b)$ is a boundary edge, then $\bold{g}_{(a,b)}=\bold{0}$. 
\end{remark}

\begin{theorem}[{\cite[Theorem 3.13]{ciliberti2024}}]\label{theorem 2}
 Let $T$ be a $\theta$-invariant triangulation of $\mathbf{P}_{2n+2}$ with oriented diameter $d$, and let $\mathcal{A}=\mathcal{A}(T)^C$ be the cluster algebra of type $C_n$ with principal coefficients in $T$. 
    Let $[a,b]$ be a $\theta$-orbit, and $x_{ab}$ the cluster variable of $\mathcal{A}$ which corresponds to $[a,b]$. Let $F_{ab}$ and $\bold{g}_{ab}$ be the $F$-polynomial and the $\bold{g}$-vector of $x_{ab}$, respectively. 
    Then $F_{ab}=F_{ab}^C$ and $\bold{g}_{ab}=\bold{g}_{ab}^C$.
\end{theorem}

\section{Modified snake graphs from $\theta$-orbits}\label{modified_sg}

In this section we associate a labeled modified snake graph $\mathcal{G}_{ab}$ to each $\theta$-orbit $[a,b]$, and prove that the perfect matching polynomial $F_{\mathcal{G}_{ab}}$ (resp. the $\bold{g}$-vector $\bold{g}_{\mathcal{G}_{ab}}$) is equal to the $F$-polynomial (resp. $\bold{g}$-vector) of the cluster variable which corresponds to $[a,b]$. The definition of $\mathcal{G}_{ab}$ has been inspired by the work of Musiker \cite{M} for type $B$ and $C$ cluster algebras.

\subsection{Type B}
\begin{definition}\label{def_mod_sg_B}
    Let $\Bar{T}=\{ \tau_1, \dots, \tau_n \}$ be a triangulation of $\mathbf{P}_{n+3}$, such that $\tau_n$ is an edge of a triangle of $\Bar{T}$ whose other two edges are boundary edges. Let $\gamma$ be a diagonal of $\mathbf{P}_{n+3}$ which is not in $\Bar{T}$. We define the \emph{labeled modified snake graph} $\hat{\mathcal{G}_{\gamma}}$ associated with $\gamma$ as the usual labeled snake graph $\mathcal{G}_{\gamma}$ of Definition \ref{def_snake_graph} with these two modifications:
    \begin{itemize}
        \item the edge with label $\tau_n$ in the tile $G_{n-1}$ is replaced by three new edges in order to have $\hat{G}_{n-1}$ homeomorphic to a hexagon in the following way:
    
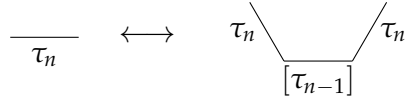
\begin{figure}[H]
        \centering
         \begin{tikzpicture}[scale=0.9]
\draw (180:1) -- node[midway, left] {$\textcolor{black}{\tau_n}$}(240:1);
\draw (240:1) --node[midway, below,yshift=0.8mm] {$\textcolor{black}{\tau_{[n-1]}}$} (300:1);
\draw (300:1) -- node[midway, right] {$\textcolor{black}{\tau_n}$}(0:1);
\begin{scope}[xshift=-4cm,yshift=0.35cm]
    \draw (240:1) -- node[midway, below] {$\textcolor{black}{\tau_n}$}(300:1);
\end{scope}
\node at (-2.5,-0.5) {$\longleftrightarrow$};
\end{tikzpicture}
    \end{figure}
    \item if $l$ is the label of an edge $e$ of $G_n$, and $e$ is an internal edge of $\mathcal{G}_\gamma$, then $l$ is also the label of the edge of $\hat{G}_n$ opposite to $e$.  
    \end{itemize}
    
\end{definition}
\begin{remark}
    In $\hat{\mathcal{G}_{\gamma}}$, unlike $\mathcal{G}_{\gamma}$, $\tau_{[n-1]}$ can also be the label of an external edge. This is the edge along which we will glue the labeled modified snake graphs of diagonals to construct the labeled modified snake graphs associated with $\theta$-orbits. See Definition \ref{lab_mod_sg_orbit}.
\end{remark}
\begin{example}\label{ex_diag_gamma}
    In the example for $n=3$ in Figure \ref{fig:snake_graph_B}, we compute snake graphs $\mathcal{G}_\gamma$ and $\hat{\mathcal{G}_\gamma}$ of a diagonal $\gamma$ in a triangulated hexagon. The tile $G_{n-1}=G_2$ is the central tile of $\mathcal{G}_\gamma$, and $\tau_n=\tau_3$ is its north edge. According to Definition \ref{def_mod_sg_B}, this edge is replaced by the three edges labeled $3,[2],3$ at the top of the central hexagon in $\hat{\mathcal{G}}_\gamma$. Furthermore, the additional label $[2]$ is added on the east edge of the last tile.   
\begin{figure}[H]
    \centering
   \begin{tikzpicture}[scale=0.9]
    \draw (0:3cm) -- (60:3cm) -- (120:3cm) -- (180:3cm) -- (240:3cm) -- (300:3cm) -- cycle;
                    \draw (0:3cm) -- node[midway, above, xshift=-1mm, yshift=-1mm] {3} (120:3cm);
                    \draw (0:3cm) -- node[midway, above left,xshift=1.2mm] {2} (180:3cm);
                    \draw (300:3cm) -- node[midway, above left,xshift=1mm] {1} (180:3cm);
                    \draw[blue, line width=0.3mm] (60:3cm) -- node[left, yshift=-9mm] {$\gamma$} (240:3cm);
   \begin{scope}[xshift=7cm]

   \draw (-0.5,0) -- (-0.5,1) -- (0.5,1) -- (0.5,0) -- cycle;

\begin{scope}
    \draw (-0.5,0) --node[midway] {$[1]$} (-0.5,1) --node[midway, above left,xshift=2mm] {3} (0.5,1) -- node[midway] {$[2]$}(0.5,0) -- node[midway, below,xshift=1mm] {1} cycle;
     \node at (0,0.5) {2};

     \node at (-2.5,0.5) {$\mathcal{G}_\gamma=$};
     \node at (-2.5,-3) {$\hat{\mathcal{G}}_{\gamma}=$};
\end{scope}

\begin{scope}[xshift=-1cm]
    \draw (-0.5,0) -- (-0.5,1) -- node[midway, above] {2}(0.5,1) -- (0.5,0) --  cycle;
    \node at (0,0.5) {1};
    
\end{scope}  

\begin{scope}[xshift=1cm]
    \draw (-0.5,0) -- (-0.5,1) --(0.5,1) -- (0.5,0) -- node[midway, below] {2}cycle;
    \node at (0,0.5) {3};
    
\end{scope} 
  \end{scope}

                      \begin{scope}[xshift=7.8cm, yshift=-2.5cm]
   \draw (0:1) -- node[midway,right] {3}(60:1) --node[midway,above,yshift=-0.8mm] {$[2]$} (120:1) -- node[midway,left] {3}(180:1);
    \draw (300:1) -- node[midway,below] {1}(240:1);

\node at (0,0) {2};

\begin{scope}[xshift=-1.5cm,yshift=-0.855cm]
    \draw[] (0:1) -- node[pos=0.5, right, xshift=-4mm, yshift=0mm] {$[1]$} (60:1);
    \draw (60:1) -- node[midway,above] {2} (120:1) -- (180:1) -- (0:1);

    \node at (0,0.4) {1};
    
\end{scope}
\begin{scope}[xshift=1.5cm,yshift=-0.855cm]
    \draw (0:1) -- node[midway] {$[2]$}(60:1) --  (120:1); 
    \draw[] (120:1) -- node[pos=0.5, left,xshift=2mm] {$[2]$}(180:1);
   \draw (180:1) -- node[midway,below] {2}(0:1);
    
    \node at (0,0.4) {3};
    
\end{scope}
\end{scope}
   \end{tikzpicture}
    \caption{The snake graphs $\mathcal{G}_\gamma$ and $\hat{\mathcal{G}_\gamma}$ for a diagonal $\gamma$ in a triangulated hexagon (type $B$).}
    \label{fig:snake_graph_B}
\end{figure}

\end{example}

\begin{remark}\label{rmk_poset_iso}
        The operation $f:Match(\mathcal{G}_\gamma)\to Match(\hat{\mathcal{G}}_\gamma)$ defined as follows:
        \begin{figure}[H]
        \centering
         \begin{tikzpicture}[scale=0.9]
\draw[red,very thick] (180:1) -- node[midway, left] {$\textcolor{black}{\tau_n}$}(240:1);
\draw (240:1) --node[midway, below,yshift=0.8mm] {$\textcolor{black}{\tau_{[n-1]}}$} (300:1);
\draw[red,very thick] (300:1) -- node[midway, right] {$\textcolor{black}{\tau_n}$}(0:1);
\begin{scope}[xshift=-4cm,yshift=0.35cm]
    \draw[red,very thick] (240:1) -- node[midway, below] {$\textcolor{black}{\tau_n}$}(300:1);
\end{scope}
\node at (-2.5,-0.5) {$\xlongleftrightarrow{f}$};
\end{tikzpicture}
    \end{figure}
       \begin{figure}[H]
        \centering
         \begin{tikzpicture}[scale=0.9]
\draw (180:1) -- node[midway, left] {$\textcolor{black}{\tau_n}$}(240:1);
\draw[red,very thick] (240:1) --node[midway, below,yshift=0.8mm] {$\textcolor{black}{\tau_{[n-1]}}$} (300:1);
\draw (300:1) -- node[midway, right] {$\textcolor{black}{\tau_n}$}(0:1);
\begin{scope}[xshift=-4cm,yshift=0.35cm]
    \draw (240:1) -- node[midway, below] {$\textcolor{black}{\tau_n}$}(300:1);
\end{scope}
\node at (-2.5,-0.5) {$\xlongleftrightarrow{f}$};
\end{tikzpicture}
    \end{figure}
    is a poset preserving isomorphism between the set of perfect matchings of $\mathcal{G}_\gamma$ and the set of perfect matchings of $\hat{\mathcal{G}}_\gamma$.
\end{remark}

\begin{definition}\label{lab_mod_sg_orbit}
    Let $T=\{ \tau_1, \dots, \tau_n=d, \dots, \tau_{2n-1} \}$ be a $\theta$-invariant triangulation of $\mathbf{P}_{2n+2}$ with oriented diameter $\tau_n=d$, such that $\tau_n$ and $\tau_{n-1}$ are edges of a triangle of $T$ whose third edge is a boundary edge. Let $[a,b]$ be a $\theta$-orbit which is not in $T$. We associate to $[a,b]$ the labeled modified snake graph $\mathcal{G}_{ab}$ defined in the following way:
    \begin{itemize}
        \item if $\text{Res}([a,b])=\{ \gamma \}$, then $\mathcal{G}_{ab}:=\hat{\mathcal{G}}_{\gamma}$;
        \item if $\text{Res}([a,b])=\{ \gamma_1,\gamma_2\}$, with $\gamma_1$ counterclockwise (resp. clockwise) from $\gamma_2$ if $\tau_{n-1}$ is counterclockwise (resp. clockwise) from $\tau_n$, then $\mathcal{G}_{ab}$ is obtained by gluing the tile with label $n$ of $\hat{\mathcal{G}}_{\gamma_2}$ to the tile with label $n-1$ of $\hat{\mathcal{G}}_{\gamma_1}$ along $\tau_{[n-1]}$, following the gluing rule recalled in Section \ref{section_sg}. If both $\hat{\mathcal{G}}_{\gamma_1}$ and $\hat{\mathcal{G}}_{\gamma_2}$ contain a tile with label $n-1$, we add an edge with label $n-1$ from the top right vertex of the tile of $\hat{\mathcal{G}}_{\gamma_1}$ with label $n$ to the top left vertex of the tile of $\hat{\mathcal{G}}_{\gamma_2}$ with label $n-1$, as in Figure \ref{additional_edge}.

    \end{itemize}

 The edges of $\mathcal{G}_{ab}$ along which two tiles are glued are called \emph{internal}, while the remaining edges are called \emph{external}.   
\end{definition}
        \begin{figure}[H]
    \centering
   \begin{tikzpicture}[scale=0.9]
   \begin{scope}[xshift=7.8cm, yshift=-2cm]
   \draw (180:1) arc (180:28:1.85) node[midway, above,xshift=3mm] {$n-1$};
   \draw (0:1) -- node[midway,right] {}(60:1) -- (120:1) -- node[midway,left] {}(180:1);
    \draw (300:1) -- node[midway,below] {}(240:1);

\node at (0,0) {$n-1$};

\begin{scope}[xshift=3cm]
\draw[] (0:1) -- node[pos=0.5, right, xshift=-4mm, yshift=0mm] {} (60:1);
    \draw (60:1) -- node[midway,above] {} (120:1) -- (180:1) --(240:1)-- (300:1) -- (0:1); 
    \node at (0,0) {$n-1$};
\end{scope}

\begin{scope}[xshift=-1.5cm,yshift=-0.855cm]
    \draw[] (0:1) -- node[pos=0.5, right, xshift=-4mm, yshift=0mm] {} (60:1);
    \draw (60:1) -- node[midway,above] {} (120:1) -- (180:1) -- (0:1);

    \node at (0,0.4) {$n$};
    
\end{scope}
\begin{scope}[xshift=1.5cm,yshift=-0.855cm]
    \draw (0:1) -- node[midway] {}(60:1) --  (120:1); 
    \draw[] (120:1) -- node[pos=0.5, left,xshift=2mm] {}(180:1);
   \draw (180:1) -- node[midway,below] {}(0:1);
    \node at (-0.3,1.6) {$I$};
    \node at (0,0.4) {$n$};
    
\end{scope}
\end{scope}
   \end{tikzpicture}
    \caption{Additional edge from the top right vertex of the tile of $\hat{\mathcal{G}}_{\gamma_1}$ with label $n$ to the top left vertex of the tile of $\hat{\mathcal{G}}_{\gamma_2}$ with label $n-1$.}
    \label{additional_edge}
\end{figure}
\begin{remark}
By \cite[Theorem 3.2]{propp2024latticestructureorientationsgraphs}, the set $L$ of perfect matchings of $\mathcal{G}_{ab}$ is a distributive lattice. We observe that $L$ is the union of two distributive lattices, the lattice of perfect matchings which contain the additional edge and the one of perfect matchings which do not contain it, connected by a single edge corresponding to the flip of the face $I$ enclosed by the additional edge. See Figure \ref{set_pm} for an example.   
\end{remark}
\begin{example}\label{ex_lab_snake_graph2_B_computation}
We compute the labeled modified snake graph $\mathcal{G}_{ab}$ of the $\theta$-orbit $[a,b]$ in Figure \ref{ex_lab_snake_graph2_}.

\begin{figure}[H]
                               \centering
                \begin{tikzpicture}[scale=0.5]
                    \draw (90:3cm) -- (135:3cm) -- (180:3cm) -- (225:3cm) -- (270:3cm) -- (315:3cm) -- (360:3cm) -- (45:3cm) -- cycle;
                    \draw (270:3cm) -- (180:3cm);;
                    \draw (90:3cm) -- (360:3cm);
                    \draw[-{Latex[length=2mm]}] (270:3cm) -- node[midway, above left,xshift=1mm] {} (90:3cm);
                    \draw (135:3cm) -- node[midway, above left,xshift=1mm] {} (270:3cm);
                    \draw (315:3cm) -- node[midway, above left,xshift=1mm] {} (90:3cm);
                    \draw[cyan, line width=0.3mm] (180:3cm) -- (45:3cm); 
                    \draw[cyan, line width=0.3mm] (225:3cm) -- (360:3cm);

                    \node at (180:3cm) [left] {$a$};
                    \node at (225:3cm) [left] {$\Bar{b}$};
                    \node at (45:3cm) [right] {$b$};
                    \node at (360:3cm) [right] {$\Bar{a}$};

                    \begin{scope}[xshift=8cm]
                      \draw (90:3cm) -- (135:3cm) -- (180:3cm) -- (225:3cm) -- (270:3cm) -- (360:3cm) -- cycle;
 
                    \draw (180:3cm) --node[midway, above left,xshift=1mm,yshift=-1mm] {1} (270:3cm);;
                    \draw (90:3cm) -- node[midway, above left,xshift=2mm,yshift=4mm] {3} (270:3cm);
                    \draw (135:3cm) -- node[midway, above left,xshift=1mm] {2} (270:3cm);

\draw[red, line width=0.3mm] (225:3cm) -- node[midway, above, xshift=-1mm, yshift=-1mm] {$\gamma_1$} (360:3cm); 
\draw[blue, line width=0.3mm] (180:3cm) -- node[midway, above, xshift=3mm, yshift=-1mm] {$\gamma_2$} (360:3cm);

\node at (180:3cm) [left] {$a$};
\node at (225:3cm) [left] {$\bar{b}$};
\node at (360:3cm) [right] {$\ast$}; 
                    \end{scope}
                  
                \end{tikzpicture}
                               \caption{A $\theta$-orbit $[a,b]$ in a triangulated octagon and its restriction (type $B_3$).}
                               \label{ex_lab_snake_graph2_}
                           \end{figure}
First, we compute $\hat{\mathcal{G}}_{\gamma_1}$ (in red) and $\hat{\mathcal{G}}_{\gamma_2}$ (in blue) from $\mathcal{G}_{\gamma_1}$ and $\mathcal{G}_{\gamma_2}$, according to Definition \ref{def_mod_sg_B}:
 \begin{figure}[H]
    \centering
   \begin{tikzpicture}[scale=0.7]
  
 
   \draw (-0.5,0) -- (-0.5,1) -- (0.5,1) -- (0.5,0) -- cycle;

\node at (-3,0.5) {$\mathcal{G}_{\gamma_1}=$};

\begin{scope}
    \draw (-0.5,0) --node[midway] {} (-0.5,1) --node[midway, above left,xshift=2mm] {} (0.5,1) -- node[midway,right,xshift=-1mm] {\textcolor{blue}{1}}(0.5,0) -- node[midway, below,yshift=1mm] {\textcolor{blue}{3}} cycle;
     \node at (0,0.5) {2};

\end{scope}

\begin{scope}[xshift=-1cm]
    \draw (-0.5,0) -- (-0.5,1) -- node[midway, above,yshift=-1mm] {\textcolor{blue}{2}}(0.5,1) -- (0.5,0) --  cycle;
    \node at (0,0.5) {3};
\end{scope}  

\begin{scope}[yshift=1cm]
    \draw (-0.5,0) --node[midway, left,xshift=1mm] {\textcolor{blue}{2}} (-0.5,1) --(0.5,1) -- (0.5,0) -- cycle;
    \node at (0,0.5) {1};
\end{scope}
 \begin{scope}[xshift=6cm,yshift=0.855cm]
 \node at (-3.5,-0.5) {$\hat{\mathcal{G}}_{\gamma_1}=$};
\begin{scope}[xshift=-0.15cm, yshift=-0.5cm ,rotate=-30]
   \filldraw[fill=red!10] (-1,0) --  (-1,1) -- (-2,1) -- (-2,0) -- node[midway,yshift=1mm] {\textcolor{blue}{2}} cycle;
     \node at (-1.5,0.5) {1}; 
\end{scope}
\begin{scope}[]
\filldraw[fill=red!10] (0:1) -- node[midway,right] {\textcolor{blue}{3}} (60:1) -- node[midway,above] {\textcolor{blue}{1}} (120:1) -- node[midway,left] {}(180:1) --(240:1)--  node[midway,below] {\textcolor{blue}{3}} (300:1) --node[midway] {\textcolor{blue}{[2]}} (0:1); 
    \node at (0,0) {$2$};
\end{scope}

\begin{scope}[xshift=-1.5cm,yshift=-0.855cm]
    \filldraw[fill=red!10] (0:1) -- node[pos=0.5, right, xshift=-4mm, yshift=0mm] {} (60:1) -- node[midway] {\textcolor{blue}{2}} (120:1) -- node[midway] {\textcolor{blue}{[2]}}(180:1) -- (0:1);
    
    \node at (0,0.4) {$3$};
    
\end{scope}

 \end{scope}
\end{tikzpicture}
\end{figure} 
 \begin{figure}[H]
    \centering
   \begin{tikzpicture}[scale=0.7]
  
 
   \draw (-0.5,0) -- (-0.5,1) -- (0.5,1) -- (0.5,0) -- cycle;

\node at (-3,0.5) {$\mathcal{G}_{\gamma_2}=$};

\begin{scope}
    \draw (-0.5,0) --node[midway] {} (-0.5,1) --node[midway, above left,xshift=2mm] {} (0.5,1) -- node[midway,right,xshift=-1mm] {\textcolor{blue}{1}}(0.5,0) -- node[midway, below,yshift=1mm] {\textcolor{blue}{3}} cycle;
     \node at (0,0.5) {2};

\end{scope}

\begin{scope}[xshift=-1cm]
    \draw (-0.5,0) -- (-0.5,1) -- node[midway, above,yshift=-1mm] {\textcolor{blue}{2}}(0.5,1) -- (0.5,0) --  cycle;
    \node at (0,0.5) {3};
\end{scope}  

 \begin{scope}[xshift=3cm,yshift=0.855cm]
 \node at (-0.5,-0.5) {$\hat{\mathcal{G}}_{\gamma_2}=$};
    \begin{scope}[xshift=3cm]
\filldraw[fill=blue!10] (0:1) -- node[midway,right] {\textcolor{blue}{3}} (60:1) -- node[midway,above] {\textcolor{blue}{1}} (120:1) -- node[midway,left] {}(180:1) --(240:1)--node[midway,below] {\textcolor{blue}{3}} (300:1) --node[midway] {\textcolor{blue}{[2]}} (0:1); 
    \node at (0,0) {$2$};
\end{scope}
\begin{scope}[xshift=1.5cm,yshift=-0.855cm]
    \filldraw[fill=blue!10] (0:1) -- node[pos=0.5, right, xshift=-4mm, yshift=0mm] {} (60:1) -- node[midway] {\textcolor{blue}{2}} (120:1) -- node[midway] {\textcolor{blue}{[2]}}(180:1) -- (0:1);
    
    \node at (0,0.4) {$3$};

\end{scope}

 \end{scope}
\end{tikzpicture}
\end{figure}
Then, according to Definition \ref{lab_mod_sg_orbit}, we glue them together and add an edge from the top right vertex of the tile of $\hat{\mathcal{G}}_{\gamma_1}$ with label 3 to the top left vertex of the tile of $\hat{\mathcal{G}}_{\gamma_2}$ with label 2. We get the following.  
\begin{figure}[H]
\centering
\begin{tikzpicture}[scale=0.9]

   \draw (-1,0) .. controls (-6,1.5) and (0,2.3) .. (2.5,0.85) node[midway, above,xshift=3mm] {$2$};

   \filldraw[fill=red!10] (0:1) -- node[midway,right] {\textcolor{blue}{3}}(60:1) -- node[midway, above] {\textcolor{blue}{1}}(120:1) -- node[midway,above] {}(180:1)--(240:1)-- node[midway,below] {\textcolor{blue}{3}}(300:1) -- (0:1);
\node at (0,0) {$2$};
\node at (-4,0) {$\mathcal{G}_{ab}=$};
\begin{scope}[xshift=-0.15cm, yshift=-0.5cm ,rotate=-30]
   \filldraw[fill=red!10] (-1,0) --  (-1,1) -- (-2,1) -- (-2,0) -- node[midway,yshift=1mm] {\textcolor{blue}{2}} cycle;
     \node at (-1.5,0.5) {1}; 
\end{scope}
\begin{scope}[xshift=3cm]
\filldraw[fill=blue!10] (0:1) -- node[midway,right] {\textcolor{blue}{3}} (60:1) -- node[midway] {\textcolor{blue}{[2]}} (120:1) -- node[midway,left] {\textcolor{blue}{3}}(180:1) --(240:1)-- (300:1) --node[midway,right] {\textcolor{blue}{1}} (0:1); 
    \node at (0,0) {$2$};
\end{scope}

\begin{scope}[xshift=-1.5cm,yshift=-0.855cm]
    \filldraw[fill=red!10] (0:1) -- node[pos=0.5, right, xshift=-4mm, yshift=0mm] {} (60:1);
    \filldraw[fill=red!10] (60:1) -- node[midway] {\textcolor{blue}{2}} (120:1) --node[midway] {\textcolor{blue}{[2]}} (180:1) -- (0:1);
    
    \node at (0,0.4) {$3$};
    
\end{scope}
\begin{scope}[xshift=1.5cm,yshift=-0.855cm]
    \filldraw[fill=blue!10] (0:1) -- node[midway] {}(60:1) --  (120:1) -- node[midway] {\textcolor{blue}{[2]}}(180:1) -- node[midway,below] {\textcolor{blue}{2}}(0:1);
    
    \node at (0,0.4) {$3$};
    
\end{scope}
\end{tikzpicture}
  \end{figure}
\end{example}
\begin{definition}
    Let $\mathcal{G}_{ab}$ be a labeled modified snake graph. We define $P_-(\mathcal{G}_{ab}) \in Match(\mathcal{G}_{ab})$ in the following way: 
    \begin{itemize}
        \item if $\mathcal{G}_{ab}=\hat{\mathcal{G}}_\gamma$, we define $P_-(\mathcal{G}_{ab}):=f(P_-(\mathcal{G}_\gamma))$, where $f$ is the bijection of Remark \ref{rmk_poset_iso};
        \item if $\mathcal{G}_{ab}$ is obtained by gluing $\hat{\mathcal{G}}_{\gamma_2}$ to $\hat{\mathcal{G}}_{\gamma_1}$, we define $P_-(\mathcal{G}_{ab}):=f(P_-(\mathcal{G}_{\gamma_1}))\cup f(P_-(\mathcal{G}_{\gamma_2}))$.
    \end{itemize}
\end{definition}
We extend the definition of height monomial $y(P)$ of a perfect matching $P$ of $\mathcal{G}_{ab}$ using Remark \ref{rmk_ind_set}.
\begin{definition}\label{def_h(P)}
Let $P_-=P_-(\mathcal{G}_{ab})$, and let $P$ be a perfect matching of $\mathcal{G}_{ab}$. The \emph{height monomial} of $P$ is
\begin{center}
    $y(P):=\displaystyle\prod_{i}y_i$,
\end{center}
where the product is over all $i$ for which $(P_- \cup P)_{|G_i}$ contains an external edge of $\mathcal{G}_{ab}$ and $P_- \cap P$ does not contain any edge of $G_i$ with label different from $\tau_n$.
\end{definition}

For a $\theta$-orbit $[a,b]$ of $\mathbf{P}_{2n+2}$ (resp. a diagonal $\gamma$ of $\mathbf{P}_{n+3}$) the minimal matching $P_-(\mathcal{G}_{ab})$ (resp. $P_-(\hat{\mathcal{G}}_{\gamma})$) is defined as in Definition \ref{p-}.

\begin{definition}\label{def_f_poly_lab_mod_sg}
    Let $[a,b]$ be a $\theta$-orbit which is not in $T$, and $\tau_{i_1}, \dots, \tau_{i_d}$ be the sequence of diagonals of $\Bar{T}=\text{Res}(T)$ crossed by the diagonals of $\text{Res}([a,b])$. Then the \emph{perfect matching polynomial of $\mathcal{G}_{ab}$} is
  \begin{center}
      $F_{\mathcal{G}_{ab}}:=\displaystyle\sum_{P}y(P)$,
  \end{center}
  where the sum is over all perfect matchings $P$ of $\mathcal{G}_{ab}$, and the $\bold{g}$-$vector$ is
  \begin{center}
      $\bold{g}_{\mathcal{G}_{ab}}:=\displaystyle\sum_{\tau_i \in P_-(\mathcal{G}_{ab})}\bold{e}_i-\displaystyle\sum_{j=1}^d \bold{e}_{i_j}$.
  \end{center}
  
The definition is extended to any $\theta$-orbit by letting $F_{\mathcal{G}_{ab}}:=1$ and $\bold{g}_{\mathcal{G}_{ab}}:=\bold{e}_i$ if $[a,b]=\{ \tau_i,\tau_{2n-i} \} \in T$, and $F_{\mathcal{G}_{ab}}:=1$ and $\bold{g}_{\mathcal{G}_{ab}}:=\bold{0}$ if $(a,b)$ is a boundary edge of $\mathbf{P}_{2n+2}$. 
\end{definition}
\begin{remark}
    In the definition of $\bold{g}_{\mathcal{G}_{ab}}$ we do not consider any external labels of $\mathcal{G}_{ab}$ of the form $\tau_{[i]}$.
\end{remark}
\begin{remark}
    The polynomial $F_{\mathcal{G}_{ab}}$ and the integer vector $\bold{g}_{\mathcal{G}_{ab}}$ uniquely define a Laurent polynomial called the \emph{perfect matching Laurent polynomial of $\mathcal{G}_{ab}$} \cite{CSI}. This Laurent polynomial can also be computed using the determinant formula presented in \cite{deloera2024determinantalformulaclustervariables}. The latter formula was important in this work, as it highlighted the need to assign a label to the additional edge of $\mathcal{G}_{ab}$.
\end{remark}
\begin{lemma}\label{lemma_bijection}
    Let $\Bar{T}=\{ \tau_1, \dots, \tau_n \}$ be a triangulation of $\mathbf{P}_{n+3}$, such that $\tau_n$ is an edge of a triangle of $\Bar{T}$ whose other two edges are boundary edges. Let $\gamma$ be a diagonal of $\mathbf{P}_{n+3}$ which is not in $\Bar{T}$. Then $F_{\hat{\mathcal{G}}_\gamma}=F_{\mathcal{G}_\gamma}$.
\end{lemma}
\begin{proof}
    Consider the bijection $f:Match(\mathcal{G}_\gamma)\to Match(\hat{\mathcal{G}}_\gamma)$ of Remark \ref{rmk_poset_iso}. We have $y(P)=y(f(P))$, for any $P \in Match(\mathcal{G}_\gamma)$. Theorefore, $F_{\hat{\mathcal{G}}_\gamma}=F_{\mathcal{G}_\gamma}$.
\end{proof}
\begin{example}
    We illustrate in Figure \ref{fig:snake_graph_1234} the proof of Lemma \ref{lemma_bijection} for the diagonal $\gamma$ of Example \ref{ex_diag_gamma}.
    \begin{figure}[H]
    \centering
   \begin{tikzpicture}[scale=0.5]
   \begin{scope}[xshift=-27cm]
        \node at (21,0.5) {$y_1y_2y_3$};
   \node at (21,-0.75) {$+$};
   \node at (21,-2) {$y_1y_3$};
   \node at (21,-3.75) {$+$};
   \node at (21,-5.5) {$y_1+y_3$};
   \node at (21,-7) {$+$};
   \node at (21,-8.5) {$1$};
   \node at (21,-10) {$\shortparallel$};
   \node at (21,-11.5) {$F_{\mathcal{G}_\gamma}$};
   \end{scope}
   \node at (21,0.5) {$y_1y_2y_3$};
   \node at (21,-0.75) {$+$};
   \node at (21,-2) {$y_1y_3$};
   \node at (21,-3.75) {$+$};
   \node at (21,-5.5) {$y_1+y_3$};
   \node at (21,-7) {$+$};
   \node at (21,-8.5) {$1$};
   \node at (21,-10) {$\shortparallel$};
   \node at (21,-11.5) {$F_{\hat{\mathcal{G}}_\gamma}$};
   \begin{scope}

\begin{scope}
    \draw (-0.5,0) -- (-0.5,1) -- (0.5,1) -- (0.5,0) --  cycle;
    \draw[red,very thick] (-0.5,1) -- (0.5,1);
   \draw[red,very thick] (0.5,0) -- (-0.5,0);
     \node at (0,0.5) {2};
     
\end{scope}

\begin{scope}[xshift=-1cm]
    \draw (-0.5,0) -- (-0.5,1) -- (0.5,1) -- (0.5,0) --  cycle;
    \draw[red,very thick] (-0.5,0) -- (-0.5,1);
    \node at (0,0.5) {1};
    
\end{scope}  

\begin{scope}[xshift=1cm]
    \draw (-0.5,0) -- (-0.5,1) -- (0.5,1) -- (0.5,0) -- cycle;
     \draw[red,very thick] (0.5,1) -- (0.5,0);
    \node at (0,0.5) {3};
    
\end{scope} 
  \end{scope}
  \draw (0,-3.5)--node[midway, above left,yshift=-1.5mm] {1}(-1,-4.5);
     \begin{scope}[xshift=-2cm, yshift=-6cm]

\begin{scope}
    \draw (-0.5,0) -- (-0.5,1) -- (0.5,1) -- (0.5,0) --  cycle;
   \draw[red,very thick] (0.5,1) -- (0.5,0);
     \node at (0,0.5) {2};
     
\end{scope}

\begin{scope}[xshift=-1cm]
    \draw (-0.5,0) -- (-0.5,1) -- (0.5,1) -- (0.5,0) --  cycle;
    \draw[red,very thick] (-0.5,1) -- (0.5,1);
    \draw[red,very thick] (0.5,0) -- (-0.5,0);
    \node at (0,0.5) {1};
    
\end{scope}  

\begin{scope}[xshift=1cm]
    \draw (-0.5,0) -- (-0.5,1) -- (0.5,1) -- (0.5,0) -- cycle;
     \draw[red,very thick] (0.5,1) -- (0.5,0);
    \node at (0,0.5) {3};
    
\end{scope} 
  \end{scope}
  \draw (0,-0.5)--node[midway, right] {2}(0,-1.5);
     \begin{scope}[yshift=-3cm]

\begin{scope}
    \draw (-0.5,0) -- (-0.5,1) -- (0.5,1) -- (0.5,0) --  cycle;
    \draw[red,very thick] (-0.5,0) -- (-0.5,1);
   \draw[red,very thick] (0.5,1) -- (0.5,0);
     \node at (0,0.5) {2};
     
\end{scope}

\begin{scope}[xshift=-1cm]
    \draw (-0.5,0) -- (-0.5,1) -- (0.5,1) -- (0.5,0) --  cycle;
    \draw[red,very thick] (-0.5,0) -- (-0.5,1);
    \node at (0,0.5) {1};
    
\end{scope}  

\begin{scope}[xshift=1cm]
    \draw (-0.5,0) -- (-0.5,1) -- (0.5,1) -- (0.5,0) -- cycle;
     \draw[red,very thick] (0.5,1) -- (0.5,0);
    \node at (0,0.5) {3};
    
\end{scope} 
  \end{scope}
    \draw (0,-3.5)--node[midway, above right, yshift=-1.5mm] {3}(1,-4.5);
     \begin{scope}[xshift=2cm, yshift=-6cm]

\begin{scope}
    \draw (-0.5,0) -- (-0.5,1) -- (0.5,1) -- (0.5,0) --  cycle;
    \draw[red,very thick] (-0.5,0) -- (-0.5,1);
     \node at (0,0.5) {2};
     
\end{scope}

\begin{scope}[xshift=-1cm]
    \draw (-0.5,0) -- (-0.5,1) -- (0.5,1) -- (0.5,0) --  cycle;
    \draw[red,very thick] (-0.5,0) -- (-0.5,1);
    \node at (0,0.5) {1};
    
\end{scope}  

\begin{scope}[xshift=1cm]
    \draw (-0.5,0) -- (-0.5,1) -- (0.5,1) -- (0.5,0) -- cycle;
     \draw[red,very thick] (-0.5,0) -- (0.5,0);
     \draw[red,very thick] (-0.5,1) -- (0.5,1);
    \node at (0,0.5) {3};
    
\end{scope} 
  \end{scope}
  \draw (0,-0.5)--node[midway, right] {2}(0,-1.5);
     \begin{scope}[yshift=-3cm]

\begin{scope}
    \draw (-0.5,0) -- (-0.5,1) -- (0.5,1) -- (0.5,0) --  cycle;
    \draw[red,very thick] (-0.5,0) -- (-0.5,1);
   \draw[red,very thick] (0.5,1) -- (0.5,0);
     \node at (0,0.5) {2};
     
\end{scope}

\begin{scope}[xshift=-1cm]
    \draw (-0.5,0) -- (-0.5,1) -- (0.5,1) -- (0.5,0) --  cycle;
    \draw[red,very thick] (-0.5,0) -- (-0.5,1);
    \node at (0,0.5) {1};
    
\end{scope}  

\begin{scope}[xshift=1cm]
    \draw (-0.5,0) -- (-0.5,1) -- (0.5,1) -- (0.5,0) -- cycle;
     \draw[red,very thick] (0.5,1) -- (0.5,0);
    \node at (0,0.5) {3};
    
\end{scope} 
  \end{scope}

\draw (0,-7.5)--node[midway, below right, yshift=1.5mm] {1}(1,-6.5);
\draw (0,-7.5)--node[midway, below left,yshift=1.5mm] {3}(-1,-6.5);

     \begin{scope}[yshift=-8.9cm]

\begin{scope}
    \draw (-0.5,0) -- (-0.5,1) -- (0.5,1) -- (0.5,0) --  cycle;
     \node at (0,0.5) {2};
     
\end{scope}

\begin{scope}[xshift=-1cm]
    \draw (-0.5,0) -- (-0.5,1) -- (0.5,1) -- (0.5,0) --  cycle;
    \draw[red,very thick] (-0.5,1) -- (0.5,1);
    \draw[red,very thick] (0.5,0) -- (-0.5,0);
    \node at (0,0.5) {1};
    
\end{scope}  

\begin{scope}[xshift=1cm]
    \draw (-0.5,0) -- (-0.5,1) -- (0.5,1) -- (0.5,0) -- cycle;
     \draw[red,very thick] (-0.5,0) -- (0.5,0);
     \draw[red,very thick] (-0.5,1) -- (0.5,1);
    \node at (0,0.5) {3};
    
\end{scope} 
  \end{scope}
  \draw (0,-0.5)--node[midway, right] {2}(0,-1.5);
     \begin{scope}[yshift=-3cm]

\begin{scope}
    \draw (-0.5,0) -- (-0.5,1) -- (0.5,1) -- (0.5,0) --  cycle;
    \draw[red,very thick] (-0.5,0) -- (-0.5,1);
   \draw[red,very thick] (0.5,1) -- (0.5,0);
     \node at (0,0.5) {2};
     
\end{scope}

\begin{scope}[xshift=-1cm]
    \draw (-0.5,0) -- (-0.5,1) -- (0.5,1) -- (0.5,0) --  cycle;
    \draw[red,very thick] (-0.5,0) -- (-0.5,1);
    \node at (0,0.5) {1};
    
\end{scope}  

\begin{scope}[xshift=1cm]
    \draw (-0.5,0) -- (-0.5,1) -- (0.5,1) -- (0.5,0) -- cycle;
     \draw[red,very thick] (0.5,1) -- (0.5,0);
    \node at (0,0.5) {3};
    
\end{scope} 
  \end{scope}

\node at (6,-3) {$\xlongleftrightarrow{f}$};
\node at (6,-11.5) {$=$};
\begin{scope}[xshift=2cm]
                        \begin{scope}[xshift=11cm,yshift=0.8cm]
                      
   \draw[red,very thick] (0:1) -- node[midway,right] {}(60:1);
   \draw (60:1)--node[midway,above,yshift=-0.8mm] {} (120:1);
   \draw[red,very thick] (120:1) -- node[midway,left] {}(180:1);
    \draw[red,very thick] (300:1) -- node[midway,below] {}(240:1);

\node at (0,0) {2};

\begin{scope}[xshift=-1.5cm,yshift=-0.855cm]
    \draw[] (0:1) -- node[pos=0.5, right, xshift=-4mm, yshift=0mm] {} (60:1);
    \draw (60:1) -- node[midway,above] {} (120:1);
    \draw[red,very thick] (120:1) -- (180:1);
    \draw (180:1)-- (0:1);

    \node at (0,0.4) {1};
    
\end{scope}
\begin{scope}[xshift=1.5cm,yshift=-0.855cm]
    \draw[red,very thick] (0:1) -- node[midway] {}(60:1);
    \draw (60:1) --  (120:1); 
    \draw[] (120:1) -- node[pos=0.5, left,xshift=2mm] {}(180:1);
   \draw (180:1) -- node[midway,below] {}(0:1);
    
    \node at (0,0.4) {3};
    
\end{scope}
\end{scope}

\draw (11,-0.4)--(11,-1.1);
\node at (11.3,-0.75) {2};
\draw (11,-3.5)--node[midway, above left,yshift=-1mm] {1}(8,-4);
\draw (11,-3.5)--node[midway, above right,yshift=-1mm] {3}(14,-4);
\draw (8,-6.3)--node[midway, left,yshift=-1mm] {3}(11,-7);
\draw (14,-6.3)--node[midway, right,yshift=-1mm] {1}(11,-7);
\begin{scope}[xshift=11cm,yshift=-2.2cm]
                      
   \draw (0:1) -- node[midway,right] {}(60:1);
   \draw[red,very thick] (60:1)--node[midway,above,yshift=-0.8mm] {} (120:1);
   \draw (120:1) -- node[midway,left] {}(180:1);
    \draw (300:1) -- node[midway,below] {}(240:1);

\node at (0,0) {2};

\begin{scope}[xshift=-1.5cm,yshift=-0.855cm]
    \draw[red,very thick] (0:1) -- node[pos=0.5, right, xshift=-4mm, yshift=0mm] {} (60:1);
    \draw (60:1) -- node[midway,above] {} (120:1);
    \draw[red,very thick] (120:1) -- (180:1);
    \draw (180:1)-- (0:1);

    \node at (0,0.4) {1};
    
\end{scope}
\begin{scope}[xshift=1.5cm,yshift=-0.855cm]
    \draw[red,very thick] (0:1) -- node[midway] {}(60:1);
    \draw (60:1) --  (120:1); 
    \draw[red,very thick] (120:1) -- node[pos=0.5, left,xshift=2mm] {}(180:1);
   \draw (180:1) -- node[midway,below] {}(0:1);
    
    \node at (0,0.4) {3};
    
\end{scope}
\end{scope}
\begin{scope}[xshift=8cm,yshift=-5.1cm]
                      
   \draw (0:1) -- node[midway,right] {}(60:1);
   \draw[red,very thick] (60:1)--node[midway,above,yshift=-0.8mm] {} (120:1);
   \draw (120:1) -- node[midway,left] {}(180:1);
    \draw (300:1) -- node[midway,below] {}(240:1);

\node at (0,0) {2};

\begin{scope}[xshift=-1.5cm,yshift=-0.855cm]
    \draw[] (0:1) -- node[pos=0.5, right, xshift=-4mm, yshift=0mm] {} (60:1);
    \draw[red,very thick] (60:1) -- node[midway,above] {} (120:1);
    \draw (120:1) -- (180:1);
    \draw[red,very thick] (180:1)-- (0:1);

    \node at (0,0.4) {1};
    
\end{scope}
\begin{scope}[xshift=1.5cm,yshift=-0.855cm]
    \draw[red,very thick] (0:1) -- node[midway] {}(60:1);
    \draw (60:1) --  (120:1); 
    \draw[red,very thick] (120:1) -- node[pos=0.5, left,xshift=2mm] {}(180:1);
   \draw (180:1) -- node[midway,below] {}(0:1);
    
    \node at (0,0.4) {3};
    
\end{scope}
\end{scope}
\begin{scope}[xshift=14cm,yshift=-5.1cm]
                      
   \draw (0:1) -- node[midway,right] {}(60:1);
   \draw[red,very thick] (60:1)--node[midway,above,yshift=-0.8mm] {} (120:1);
   \draw (120:1) -- node[midway,left] {}(180:1);
    \draw (300:1) -- node[midway,below] {}(240:1);

\node at (0,0) {2};

\begin{scope}[xshift=-1.5cm,yshift=-0.855cm]
    \draw[red,very thick] (0:1) -- node[pos=0.5, right, xshift=-4mm, yshift=0mm] {} (60:1);
    \draw (60:1) -- node[midway,above] {} (120:1);
    \draw[red,very thick] (120:1) -- (180:1);
    \draw (180:1)-- (0:1);

    \node at (0,0.4) {1};
    
\end{scope}
\begin{scope}[xshift=1.5cm,yshift=-0.855cm]
    \draw (0:1) -- node[midway] {}(60:1);
    \draw[red,very thick] (60:1) --  (120:1); 
    \draw[] (120:1) -- node[pos=0.5, left,xshift=2mm] {}(180:1);
   \draw[red,very thick] (180:1) -- node[midway,below] {}(0:1);
    
    \node at (0,0.4) {3};
    
\end{scope}
\end{scope}
\begin{scope}[xshift=11cm,yshift=-8cm]
                      
   \draw (0:1) -- node[midway,right] {}(60:1);
   \draw[red,very thick] (60:1)--node[midway,above,yshift=-0.8mm] {} (120:1);
   \draw (120:1) -- node[midway,left] {}(180:1);
    \draw (300:1) -- node[midway,below] {}(240:1);

\node at (0,0) {2};

\begin{scope}[xshift=-1.5cm,yshift=-0.855cm]
    \draw[] (0:1) -- node[pos=0.5, right, xshift=-4mm, yshift=0mm] {} (60:1);
    \draw[red,very thick] (60:1) -- node[midway,above] {} (120:1);
    \draw (120:1) -- (180:1);
    \draw[red,very thick] (180:1)-- (0:1);

    \node at (0,0.4) {1};
    
\end{scope}
\begin{scope}[xshift=1.5cm,yshift=-0.855cm]
    \draw (0:1) -- node[midway] {}(60:1);
    \draw[red,very thick] (60:1) --  (120:1); 
    \draw[] (120:1) -- node[pos=0.5, left,xshift=2mm] {}(180:1);
   \draw[red,very thick] (180:1) -- node[midway,below] {}(0:1);
    
    \node at (0,0.4) {3};
    
\end{scope}
\end{scope}  
\end{scope}    

   \end{tikzpicture}
    \caption{The posets of perfect matchings of $\mathcal{G}_\gamma$ and $\hat{\mathcal{G}_\gamma}$, and the corresponding monomials which give $F_{\mathcal{G}_\gamma}=F_{\hat{\mathcal{G}}_\gamma}$.}
    \label{fig:snake_graph_1234}
\end{figure}
\end{example}
\begin{example}\label{ex_lab_snake_graph2_B}
We compute the perfect matching polynomial $F_{\mathcal{G}_{ab}}$ and the $\bold{g}$-vector $\bold{g}_{\mathcal{G}_{ab}}$ of the labeled modified snake graph $\mathcal{G}_{ab}$ of Example \ref{ex_lab_snake_graph2_B_computation}. The poset of all perfect matchings of $\mathcal{G}_{ab}$, with the corresponding monomials, is illustrated in Figure \ref{set_pm}.
 
We have 
\begin{align*}
    F_{\mathcal{G}_{ab}}=y_1y_2^2y_3^2+y_2^2y_3^2+y_1y_2y_3^2+2y_2y_3^2+2y_2y_3+y_3^2+2y_3+1,
\end{align*}
and 
\begin{align*}
    \bold{g}_{\mathcal{G}_{ab}}=\begin{pmatrix} 2\\ 2\\ 0
    \end{pmatrix}-\begin{pmatrix} 1\\ 2\\ 2
    \end{pmatrix}=\begin{pmatrix} 1\\ 0\\ -2
    \end{pmatrix}.
\end{align*}
\end{example}
\begin{lemma}\label{lemma1_B}
   Let $T=\{ \tau_1, \dots, \tau_n=d, \dots, \tau_{2n-1} \}$ be a $\theta$-invariant triangulation of $\mathbf{P}_{2n+2}$ such that $\tau_n$ and $\tau_{n-1}$ are edges of a triangle of $T$ whose third edge is a boundary edge, and $\tau_n=d$ is oriented. For any $\theta$-orbit $[a,b]$ of $\mathbf{P}_{2n+2}$, $F_{\mathcal{G}_{ab}}=F^B_{ab}$ (cf. Definition \ref{def_type_B}). 
\end{lemma}

\begin{proof}
    If $\text{Res}([a,b])=\{\gamma\}$, the statement holds since $F_{\mathcal{G}_{ab}}=F_{\hat{\mathcal{G}}_{\gamma}}=F_{\mathcal{G}_{\gamma}}$ (cf. Lemma \ref{lemma_bijection}). Otherwise, if $\text{Res}([a,b])=\{\gamma_1,\gamma_2\}$, we have two cases to consider.
    \begin{itemize}
        \item [i)] One of $\gamma_1$ and $\gamma_2$, say $\gamma_2$, intersects only $\tau_n$. So $\mathcal{G}_{ab}$ will be of the following form:
             \begin{figure}[H]
    \centering
   \begin{tikzpicture}[scale=1]
   \begin{scope}[xshift=7.8cm, yshift=-2cm]
   \filldraw[fill=red!10] (0:1) -- node[midway,right] {}(60:1) --node[midway, above] {$\vdots$} (120:1) -- node[midway, above left] {$\ddots$}(180:1) --  (180:1) -- (240:1) -- (300:1);

\node at (0,0) {$n-1$};

\begin{scope}[xshift=-1.5cm,yshift=-0.855cm]
    \filldraw[fill=red!10] (0:1) -- node[pos=0.5, right, xshift=-4mm, yshift=0mm] {} (60:1);
    \filldraw[fill=red!10] (60:1) -- node[midway,above] {} (120:1) -- (180:1) -- (0:1);

    \node at (0,0.4) {$n$};
    
\end{scope}
\begin{scope}[xshift=1.5cm,yshift=-0.855cm]
   \filldraw[fill=blue!10] (0:1) -- node[pos=0.5, right, xshift=-4mm, yshift=0mm] {} (60:1);
    \filldraw[fill=blue!10] (60:1) -- node[midway,above] {} (120:1) -- node[pos=0.5, right, xshift=-6mm, yshift=-2mm] {$\tau_{[n-1]}$}  (180:1) -- (0:1);
    
    \node at (0,0.4) {$n$};
    
\end{scope}
\end{scope}
   \end{tikzpicture}
    \label{}
\end{figure} 
where the red part represents $\hat{\mathcal{G}}_{\gamma_1}$, and the blue part represents $\hat{\mathcal{G}}_{\gamma_2}$. We have that
\begin{equation}
    F_{\mathcal{G}_{ab}}=F_{\hat{\mathcal{G}}_{\gamma_1}}F_{\hat{\mathcal{G}}_{\gamma_2}}-R=F_{\mathcal{G}_{\gamma_1}}F_{\mathcal{G}_{\gamma_2}}-R,
\end{equation}
where $R$ is the sum of the monomials which correspond to gluing of perfect matchings of $\hat{\mathcal{G}}_{\gamma_1}$ and perfect matchings of $\hat{\mathcal{G}}_{\gamma_2}$ which are not perfect matchings of $\mathcal{G}_{ab}$. They are all of the form 
             \begin{figure}[H]
    \centering
   \begin{tikzpicture}[scale=1]
   \begin{scope}[xshift=7.8cm, yshift=-2cm]
   \filldraw[fill=red!10] (0:1) -- node[midway,right] {}(60:1) --node[midway, above] {$\vdots$} (120:1) -- node[midway, above left] {$\ddots$}(180:1) -- (240:1) -- (300:1);
\draw[red,very thick] (0:1) -- node[midway, right] {$\textcolor{black}{e_1}$}(60:1);
\draw[red,very thick] (120:1) -- node[midway,left] {$\textcolor{black}{e_2}$}(180:1);
\draw[red,very thick] (240:1) -- node[midway,below] {$\textcolor{black}{e_3}$}(300:1);
\node at (0,0) {$n-1$};

\begin{scope}[xshift=-1.5cm,yshift=-0.855cm]
    \filldraw[fill=red!10] (0:1) -- node[pos=0.5, right, xshift=-4mm, yshift=0mm] {} (60:1);
    \filldraw[fill=red!10] (60:1) -- node[midway,above] {} (120:1) -- (180:1) -- (0:1);
\draw[red,very thick] (120:1) -- node[midway] {}(180:1);
    
    \node at (0,0.4) {$n$};
    
\end{scope}
\begin{scope}[xshift=1.5cm,yshift=-0.855cm]
    \filldraw[fill=blue!10] (0:1) -- node[midway] {}(60:1) --  (120:1) -- node[pos=0.5, right, xshift=-6mm, yshift=-2mm] {$\tau_{[n-1]}$}(180:1);
   \draw (180:1) -- node[midway,below] {}(0:1);
    \draw[red,very thick] (120:1) -- node[midway] {}(60:1);
    \draw[red,very thick] (180:1) -- node[midway] {}(0:1);
    \node at (0,0.4) {$n$};
    
\end{scope}
\end{scope}
   \end{tikzpicture}
    \label{forbidden_2}
\end{figure} 
Therefore, we have to describe all perfect matchings of $\hat{\mathcal{G}}_{\gamma_1}$ which do not contain the edge with label $\tau_{[n-1]}$, along which we glue $\hat{G}_n$ of $\hat{\mathcal{G}}_{\gamma_2}$ and $\hat{G}_{n-1}$ of $\hat{\mathcal{G}}_{\gamma_1}$. A similar question for $\hat{\mathcal{G}}_{\gamma_2}$ is trivial, as it is just one tile.

We consider the type $A$ exchange relation corresponding to the crossing of diagonals $\tau_{n-1}$ and $\tilde\gamma_1$, which is the diagonal of $\mathbf{P}_{n+3}$ which intersects the same diagonals of $\Bar{T}$ as $\gamma_1$ but $\tau_n$. We have two cases to consider.
\begin{itemize}
    \item [1)] $\tau_{[n-1]}$ is not in $P_-(\hat{\mathcal{G}}_{\gamma_1})$. So the red edges $e_1$, $e_2$, $e_3$ of $\hat{G}_{n-1}$ in Figure \ref{forbidden_2} are in $P_-(\hat{\mathcal{G}}_{\gamma_1})$. It follows that $y_{n-1}$ is a summand of $F_{\mathcal{G}_{\gamma_1}}$. Therefore, $\tau_{n-1}$ has to be counterclockwise from $\tau_n$.
    \begin{figure}[H]
                               \centering
                \begin{tikzpicture}[scale=0.5]
                        \draw (90:3cm) -- (110:3cm) -- (130:3cm) -- (150:3cm)-- (170:3cm) -- (190:3cm) -- (210:3cm)--(230:3cm)--(250:3cm) -- (270:3cm) -- (360:3cm)  --  cycle;
                    \draw (90:3cm) -- node[midway, right, xshift=-1mm, yshift=-1mm] {$\tau_n$} (270:3cm);
                    \draw (90:3cm) -- node[midway, above left,xshift=1mm] {} (270:3cm);
                   
                    \draw[cyan, line width=0.3mm] (90:3cm) -- node[midway, xshift=-2mm] {$\tau_{n-1}$}(250:3cm);
                     \draw[red, line width=0.3mm] (150:3cm) -- node[midway, right,yshift=-2mm] {$\tilde\gamma_1$}(270:3cm);

                    \node at (90:3cm) [above] {$c$};
                    \node at (270:3cm) [below] {$\Bar{c}$};

\node at (150:3cm) [left] {$a$};
\node at (240:3cm) [below] {$\Bar{b}$};
                    \node at (360:3cm) [right] {$\ast$};
                      
                \end{tikzpicture}
                               \caption{Type $A$ exchange relation corresponding to the crossing of $\tilde\gamma_1$ and $\tau_{n-1}$.}
                               \label{fig:enter-label}
                           \end{figure} 
                           We have that
\begin{equation}\label{eq_3_}
    F_{\tilde\gamma_1}=\bold{y}^{\bold{d}_{a\Bar{b},c\Bar{c}}}F_{(a,c)}+F_{(a,\Bar{b})}.
\end{equation}
Since $e_1$, $e_2$, $e_3$ are in the minimal perfect matching of $\hat{\mathcal{G}}_{\gamma_1}$, and so of $\hat{\mathcal{G}}_{\tilde\gamma_1}$, the sum of the monomials which correspond to the perfect matchings of $\hat{\mathcal{G}}_{\tilde\gamma_1}$ which contain them in the right hand side of \ref{eq_3_} is $F_{(a,\Bar{b})}$. 
Therefore,
\begin{equation}
    R=y_nF_{(a,\Bar{b})}=\bold{y}^{\bold{d}_{a\ast,\Bar{b}\ast}}F_{(a,\Bar{b})}.
\end{equation}

So, we obtain 
\begin{equation*}
   F_{\mathcal{G}_{ab}}=F_{\gamma_1}F_{\gamma_2}-\bold{y}^{\bold{d}_{\gamma_1,\gamma_2}}F_{(a,\Bar{b})}=F_{ab}^B. 
\end{equation*}
\item [2)] $\tau_{[n-1]}$ is in $P_-(\hat{\mathcal{G}}_{\gamma_1})$. So the edges $e_1$, $e_2$, $e_3$ in Figure \ref{forbidden_2} are not in $P_-(\hat{\mathcal{G}}_{\gamma_1})$. It follows that $y_{n-1}$ is not a summand of $F_{\mathcal{G}_{\gamma_1}}$. Therefore, $\tau_{n-1}$ has to be clockwise from $\tau_n$.

\begin{figure}[H]
                               \centering
                \begin{tikzpicture}[scale=0.5]
                        \draw (90:3cm) -- (110:3cm) -- (130:3cm) -- (150:3cm)-- (170:3cm) -- (190:3cm) -- (210:3cm)--(230:3cm)--(250:3cm) -- (270:3cm) -- (360:3cm)  --  cycle;
                    \draw (90:3cm) -- node[midway, above left,xshift=1mm] {} (270:3cm);
                    \draw (90:3cm) -- node[midway, above right, xshift=-1mm] {$\tau_n$} (270:3cm);
                   
                    \draw[cyan, line width=0.3mm] (110:3cm) -- node[midway, xshift=-2mm] {$\tau_{n-1}$}(270:3cm);
                     \draw[red, line width=0.3mm] (150:3cm) -- node[midway, below] {$\tilde\gamma_1$}(90:3cm);

                    \node at (90:3cm) [above] {$c$};
                    \node at (270:3cm) [below] {$\Bar{c}$};

\node at (150:3cm) [left] {$a$};
\node at (120:3cm) [above] {$\Bar{b}$};
                    \node at (360:3cm) [right] {$\ast$};
                      
                \end{tikzpicture}
                               \caption{Type $A$ exchange relation corresponding to the crossing of $\tilde\gamma_1$ and $\tau_{n-1}$.}
                               \label{fig:enter-label}
                           \end{figure} 
   
                           We have that
\begin{equation}\label{eq_4}
    F_{\tilde\gamma_1}=\bold{y}^{\bold{d}_{a\Bar{c},\Bar{b}c}}F_{(a,\Bar{b})}+F_{(a,\Bar{c})}.
\end{equation}
Since $e_1$, $e_2$, $e_3$ are not in the minimal perfect matching of $\hat{\mathcal{G}}_{\gamma_1}$, and so of $\hat{\mathcal{G}}_{\tilde\gamma_1}$, the sum of the monomials which correspond to the perfect matchings of $\hat{\mathcal{G}}_{\tilde\gamma_1}$ which contain them in the right hand side of \ref{eq_4} is $\bold{y}^{\bold{d}_{a\Bar{c},\Bar{b}c}}F_{(a,\Bar{b})}$. 
Therefore,
\begin{equation}
    R=y_n\bold{y}^{\bold{d}_{a\Bar{c},\Bar{b}c}}F_{(a,\Bar{b})}=\bold{y}^{\bold{d}_{a\ast,\Bar{b}\ast}}F_{(a,\Bar{b})}.
\end{equation}

So, we obtain 
\begin{equation*}
   F_{\mathcal{G}_{ab}}=F_{\gamma_1}F_{\gamma_2}-\bold{y}^{\bold{d}_{\gamma_1,\gamma_2}}F_{(a,\Bar{b})}=F_{ab}^B. 
\end{equation*}
\end{itemize}
\item [ii)] Both $\gamma_1$ and $\gamma_2$ intersect $\tau_{n-1}$. So $\mathcal{G}_{ab}$ will be of the following form:
 \begin{figure}[H]
    \centering
   \begin{tikzpicture}[scale=1]
   \begin{scope}[xshift=7.8cm, yshift=-2cm]
   \draw (180:1) arc (180:28:1.85);

    \filldraw[fill=red!10] (0:1) -- node[midway,right] {}(60:1) --node[midway, above] {$\vdots$} (120:1) -- node[midway, above left] {$\ddots$}(180:1) -- (240:1) -- (300:1);

\node at (0,0) {$n-1$};

\begin{scope}[xshift=3cm]
\filldraw[fill=blue!10] (0:1) -- node[pos=0.5, right, xshift=-4mm, yshift=0mm] {} (60:1) -- node[midway,above] {} (120:1) -- (180:1) --(240:1)-- (300:1) -- (0:1); 
     \draw (0:1) -- node[midway,right] {}(60:1) --node[midway, above] {$\vdots$} (120:1) -- node[midway, above left] {$\ddots$}(180:1);
    \node at (0,0) {$n-1$};
\end{scope}

\begin{scope}[xshift=-1.5cm,yshift=-0.855cm]
    \filldraw[fill=red!10] (0:1) --(60:1) --  (120:1) -- (180:1) -- (0:1);

    \node at (0,0.4) {$n$};
    
\end{scope}
\begin{scope}[xshift=1.5cm,yshift=-0.855cm]
    \filldraw[fill=blue!10] (0:1) -- node[midway] {}(60:1) --  (120:1)-- node[pos=0.5, right, xshift=-6mm, yshift=-2mm] {$\tau_{[n-1]}$}(180:1) -- node[midway,below] {}(0:1);
    
    \node at (0,0.4) {$n$};
    
\end{scope}
\end{scope}
   \end{tikzpicture}
    \label{additional_edge1}
\end{figure}
where the red part represents $\hat{\mathcal{G}}_{\gamma_1}$, and the blue part represents $\hat{\mathcal{G}}_{\gamma_2}$. We have that
\begin{equation}
    F_{\mathcal{G}_{ab}}=F_{\hat{\mathcal{G}}_{\gamma_1}}F_{\hat{\mathcal{G}}_{\gamma_2}}-R + S=F_{\mathcal{G}_{\gamma_1}}F_{\mathcal{G}_{\gamma_2}}-R+S,
\end{equation}
where $R$ is the sum of the monomials which correspond to gluing of perfect matchings of $\hat{\mathcal{G}}_{\gamma_1}$ and perfect matchings of $\hat{\mathcal{G}}_{\gamma_2}$ which are not perfect matchings of $\mathcal{G}_{ab}$, and so they are of the form 

 \begin{figure}[H]
    \centering
   \begin{tikzpicture}[scale=1]
   \begin{scope}[xshift=7.8cm, yshift=-2cm]
   \draw (180:1) arc (180:28:1.85);
    \filldraw[fill=red!10] (0:1) -- node[midway,right] {}(60:1) --node[midway, above] {$\vdots$} (120:1) -- node[midway, above left] {$\ddots$}(180:1)-- (240:1) -- (300:1);
\draw[red,very thick] (0:1) -- node[right] {\textcolor{black}{$e_1$}}(60:1);
\draw[red,very thick] (120:1) -- node[left] {\textcolor{black}{$e_2$}}(180:1);
\draw[red,very thick] (240:1) -- node[below] {\textcolor{black}{$e_3$}}(300:1);
\node at (0,0) {$n-1$};

\begin{scope}[xshift=3cm]
\filldraw[fill=blue!10] (0:1) -- node[pos=0.5, right, xshift=-4mm, yshift=0mm] {}  (60:1) -- node[midway,above] {} (120:1) -- (180:1) --(240:1)-- (300:1) -- (0:1); 
     \draw (0:1) -- node[midway,right] {}(60:1) --node[midway, above] {$\vdots$} (120:1) -- node[midway, above left] {$\ddots$}(180:1);
    \node at (0,0) {$n-1$};
    \draw[red,very thick] (120:1) -- node[below] {\textcolor{black}{$e_4$}}(60:1);
\draw[red,very thick] (0:1) -- node[right] {\textcolor{black}{$e_5$}}(300:1);
\end{scope}

\begin{scope}[xshift=-1.5cm,yshift=-0.855cm]
    \filldraw[fill=red!10] (0:1) -- node[pos=0.5, right, xshift=-4mm, yshift=0mm] {} (60:1)-- node[midway,above] {} (120:1) -- (180:1) -- (0:1);

 \draw[red,very thick] (120:1) -- node[midway] {}(180:1);   
    \node at (0,0.4) {$n$};
    
\end{scope}
\begin{scope}[xshift=1.5cm,yshift=-0.855cm]
    \filldraw[fill=blue!10] (0:1) -- node[midway] {}(60:1) --  (120:1) -- node[pos=0.5, right, xshift=-6mm, yshift=-2mm] {$\tau_{[n-1]}$}(180:1)-- node[midway,below] {}(0:1);
  \draw[red,very thick] (0:1) -- node[midway] {}(180:1);
   \draw[red,very thick] (120:1) -- node[midway] {}(60:1);
    \node at (0,0.4) {$n$};
    
\end{scope}
\end{scope}
   \end{tikzpicture}
    \label{additional_edge_A}
\end{figure}

while $S$ is the sum of the monomials which correspond to perfect matchings of $\mathcal{G}_{ab}$ which contain the additional edge, and so they are of the form
 \begin{figure}[H]
    \centering
   \begin{tikzpicture}[scale=1]
   \begin{scope}[xshift=7.8cm, yshift=-2cm]
   \draw[red,very thick] (180:1) arc (180:28:1.85);

    \filldraw[fill=red!10] (0:1) -- node[midway,right] {}(60:1) --node[midway, above] {$\vdots$} (120:1) -- node[midway, above left] {$\ddots$}(180:1) -- (240:1) -- (300:1);
\draw[red,very thick] (120:1) -- node[midway] {}(60:1);
\draw[red,very thick] (240:1) -- node[midway] {}(300:1);
\node at (0,0) {$n-1$};
\begin{scope}[xshift=3cm]
\filldraw[fill=blue!10] (0:1) -- node[pos=0.5, right, xshift=-4mm, yshift=0mm] {} (60:1) -- node[midway,above] {} (120:1) -- (180:1) --(240:1)-- (300:1) -- (0:1); 
     \draw (0:1) -- node[midway,right] {}(60:1) --node[midway, above] {$\vdots$} (120:1) -- node[midway, above left] {$\ddots$}(180:1);
    \node at (0,0) {$n-1$};
    \draw[red,very thick] (0:1) -- node[midway] {}(60:1);
\draw[red,very thick] (240:1) -- node[midway] {}(300:1);
    
\end{scope}

\begin{scope}[xshift=-1.5cm,yshift=-0.855cm]
    \filldraw[fill=red!10] (0:1) -- node[pos=0.5, right, xshift=-4mm, yshift=0mm] {} (60:1) -- node[midway,above] {} (120:1) -- (180:1) -- (0:1);

 \draw[red,very thick] (120:1) -- node[midway] {}(180:1);   
    \node at (0,0.4) {$n$};
    
\end{scope}
\begin{scope}[xshift=1.5cm,yshift=-0.855cm]
    \filldraw[fill=blue!10] (0:1) -- node[midway] {}(60:1) --  (120:1) -- node[pos=0.5, right, xshift=-6mm, yshift=-2mm] {$\tau_{[n-1]}$}(180:1)-- node[midway,below] {}(0:1);
   \draw[red,very thick] (120:1) -- node[midway] {}(60:1);
    \node at (0,0.4) {$n$};
    
\end{scope}
\end{scope}
   \end{tikzpicture}
    \label{additional_edge_B}
\end{figure}

Therefore, first we have to describe all perfect matchings of $\hat{\mathcal{G}}_{\gamma_1}$ and $\hat{\mathcal{G}}_{\gamma_2}$ which do not contain the edge with label $\tau_{[n-1]}$, along which we glue $\hat{G}_n$ of $\hat{\mathcal{G}}_{\gamma_2}$ and $\hat{G}_{n-1}$ of $\hat{\mathcal{G}}_{\gamma_1}$. 

We consider the type $A$ exchange relation corresponding to the crossing of diagonals $\tau_{n-1}$ and $\tilde\gamma_1$ (resp. $\tilde\gamma_2$), which is the diagonal of $\mathbf{P}_{n+3}$ which intersects the same diagonals of $\Bar{T}$ as $\gamma_1$ (resp. $\gamma_2$) but $\tau_n$. We have two cases to consider.
\begin{itemize}
    \item [1)] $\tau_{[n-1]}$ is not in $P_-(\hat{\mathcal{G}}_{\gamma_1})$. So the red edges $e_1$, $e_2$, $e_3$ of the tile of $\hat{\mathcal{G}}_{\gamma_1}$ with label $n-1$ in Figure \ref{additional_edge_A} are in $P_-(\hat{\mathcal{G}}_{\gamma_1})$. It follows that $y_{n-1}$ is a summand of $F_{\mathcal{G}_{\gamma_1}}$. Therefore, $\tau_{n-1}$ has to be counterclockwise from $\tau_n$.
    \begin{figure}[H]
                               \centering
                \begin{tikzpicture}[scale=0.5]
                       \draw (90:3cm) -- (110:3cm) -- (130:3cm) -- (150:3cm)-- (170:3cm) -- (190:3cm) -- (210:3cm)--(230:3cm)--(250:3cm) -- (270:3cm) -- (360:3cm)  --  cycle;
                    \draw (90:3cm) -- node[midway, right, xshift=-1mm, yshift=-1mm] {$\tau_n$} (270:3cm);
                    \draw (90:3cm) -- node[midway, above left,xshift=1mm] {} (270:3cm);
                   
                    \draw[cyan, line width=0.3mm] (90:3cm) -- node[midway, xshift=-2mm] {$\tau_{n-1}$}(250:3cm);
                     \draw[red, line width=0.3mm] (150:3cm) -- node[midway, right,yshift=-2mm] {$\tilde\gamma_1$}(270:3cm);

                    \node at (90:3cm) [above] {$c$};
                    \node at (270:3cm) [below] {$\Bar{c}$};

\node at (150:3cm) [left] {$a$};
\node at (210:3cm) [left] {$\Bar{b}$};
\node at (240:3cm) [below] {$d$};
                    \node at (360:3cm) [right] {$\ast$};
                      
                \end{tikzpicture}
                               \caption{Type $A$ exchange relation corresponding to the crossing of $\tilde\gamma_1$ and $\tau_{n-1}$.}
                               \label{fig:enter-label}
                           \end{figure}

    \begin{figure}[H]
                               \centering
                \begin{tikzpicture}[scale=0.5]
                        \draw (90:3cm) -- (110:3cm) -- (130:3cm) -- (150:3cm)-- (170:3cm) -- (190:3cm) -- (210:3cm)--(230:3cm)--(250:3cm) -- (270:3cm) -- (360:3cm)  --  cycle;
                    \draw (90:3cm) -- node[midway, right, xshift=-1mm, yshift=-1mm] {$\tau_n$} (270:3cm);
                    \draw (90:3cm) -- node[midway, above left,xshift=1mm] {} (270:3cm);
                   
                    \draw[cyan, line width=0.3mm] (90:3cm) -- node[midway, xshift=-2mm] {$\tau_{n-1}$}(250:3cm);
                     \draw[blue, line width=0.3mm] (210:3cm) -- node[midway, above] {$\tilde\gamma_2$}(270:3cm);

                    \node at (90:3cm) [above] {$c$};
                    \node at (270:3cm) [below] {$\Bar{c}$};

\node at (150:3cm) [left] {$a$};
\node at (210:3cm) [left] {$\Bar{b}$};
\node at (240:3cm) [below] {$d$};
                    \node at (360:3cm) [right] {$\ast$};
                      
                \end{tikzpicture}
                               \caption{Type $A$ exchange relation corresponding to the crossing of $\tilde\gamma_2$ and $\tau_{n-1}$.}
                               \label{fig:enter-label}
                           \end{figure}
                           We have 
\begin{equation}\label{eq_3}
    F_{\tilde\gamma_1}=\bold{y}^{\bold{d}_{ad,c\Bar{c}}}F_{(a,c)}+F_{(a,d)},
\end{equation}
                           and
\begin{equation}\label{eq_5}
    F_{\tilde\gamma_2}=\bold{y}^{\bold{d}_{\Bar{b}d,c\Bar{c}}}F_{(\Bar{b},c)}+F_{(\Bar{b},d)}.
\end{equation}

Since $e_1$, $e_2$, $e_3$ are in $P_-(\hat{\mathcal{G}}_{\gamma_1})$, and so in $P_-(\hat{\mathcal{G}}_{\tilde\gamma_1})$, the sum of the monomials which correspond to the perfect matchings of $\hat{\mathcal{G}}_{\tilde\gamma_1}$ which contain them in the right hand side of \ref{eq_3} is $F_{(a,d)}$. 
Since $e_1$, $e_2$, $e_3$ are in $P_-(\hat{\mathcal{G}}_{\gamma_1})$, it also follows that the red edges $e_4$, $e_5$ of the tile of $\hat{\mathcal{G}}_{\gamma_2}$ with label $n-1$ in Figure \ref{additional_edge_A} are not in $P_-(\hat{\mathcal{G}}_{\gamma_2})$. So the sum of the monomials which correspond to the perfect matchings of $\hat{\mathcal{G}}_{\tilde\gamma_2}$ which contain $e_4$, $e_5$ in the right hand side of \ref{eq_5} is $\bold{y}^{\bold{d}_{\Bar{b}d,c\Bar{c}}}F_{(\Bar{b},c)}$. 
Therefore,
\begin{equation}
    R=y_n\bold{y}^{\bold{d}_{\Bar{b}d,c\Bar{c}}}F_{(a,d)}F_{(\Bar{b},c)}.
\end{equation}
On the other hand, let $f_1$, $f_2$ be the red edges of the tile of $\hat{\mathcal{G}}_{\gamma_1}$ with label $n-1$ in Figure \ref{additional_edge_B}. Then $f_2=e_3$ has label $\tau_n$, while $f_1$ has a different label. Since $e_1$, $e_2$, $e_3$ are in $P_-(\hat{\mathcal{G}}_{\gamma_1})$, it follows that $f_1$ is not in $P_-(\hat{\mathcal{G}}_{\gamma_1})$, and so it is not in $P_-(\hat{\mathcal{G}}_{\tilde\gamma_1})$. So, if $P$ is a perfect matching of $\hat{\mathcal{G}}_{\tilde\gamma_1}$ which contain $f_1$, then $h(P)$ is a multiple of $y_{n-1}$ (cf. Definition \ref{def_h(P)}). Therefore, the sum of the monomials which correspond to the perfect matchings of $\hat{\mathcal{G}}_{\tilde\gamma_1}$ which contain $f_1$ in the right hand side of \ref{eq_3} is $\bold{y}^{\bold{d}_{ad,c\Bar{c}}}F_{(a,c)}$. Moreover, since $e_1$, $e_2$, $e_3$ are in $P_-(\hat{\mathcal{G}}_{\gamma_1})$, it also follows that $f_3$, $f_4$ are in $P_-(\hat{\mathcal{G}}_{\gamma_2})$, and so in $P_-(\hat{\mathcal{G}}_{\tilde\gamma_2})$. Then the sum of the monomials which correspond to the perfect matchings of $\hat{\mathcal{G}}_{\tilde\gamma_2}$ which contain $f_3$, $f_4$ in the right hand side of \ref{eq_5} is $F_{(\Bar{b},d)}$. Therefore, 
\begin{equation}
S=y_n\bold{y}^{\bold{d}_{ad,c\Bar{c}}}F_{(a,c)}F_{(\Bar{b},d)}.
\end{equation}
Finally, we consider the exchange relation corresponding to the crossing of $(a,d)$ and $(\Bar{b},c)$.

 \begin{figure}[H]
                               \centering
                \begin{tikzpicture}[scale=0.5]
                       \draw (90:3cm) -- (110:3cm) -- (130:3cm) -- (150:3cm)-- (170:3cm) -- (190:3cm) -- (210:3cm)--(230:3cm)--(250:3cm) -- (270:3cm) -- (360:3cm)  --  cycle;
                    \draw (90:3cm) -- node[midway, right, xshift=-1mm, yshift=-1mm] {$\tau_n$} (270:3cm);
                    \draw (90:3cm) -- node[midway, above left,xshift=1mm] {} (270:3cm);
                   
                    \draw[green, line width=0.3mm] (90:3cm) -- node[midway, xshift=-2mm] {}(210:3cm);
                    \draw[cyan, line width=0.3mm] (90:3cm) -- node[midway, xshift=-2mm] {$\tau_{n-1}$}(250:3cm);
                     \draw[yellow, line width=0.3mm] (150:3cm) -- node[midway, right,yshift=-2mm] {}(250:3cm);

                    \node at (90:3cm) [above] {$c$};
                    \node at (270:3cm) [below] {$\Bar{c}$};

\node at (150:3cm) [left] {$a$};
\node at (210:3cm) [left] {$\Bar{b}$};
\node at (240:3cm) [below] {$d$};
                    \node at (360:3cm) [right] {$\ast$};
                      
                \end{tikzpicture}
                               \caption{Type $A$ exchange relation corresponding to the crossing of $(a,d)$ and $(\Bar{b},c)$.}
                               \label{fig:enter-label}
                           \end{figure} 
We have that
\begin{equation}
    F_{(a,d)}F_{(\Bar{b},c)}=\bold{y}^{\bold{d}_{ac,\Bar{b}d}}F_{(a,\Bar{b})}+ \bold{y}^{\bold{d}_{a\Bar{b},cd}}F_{(a,c)}F_{(\Bar{b},d)}.
\end{equation}

Therefore,
\begin{equation}
    -R+S=-y_n\bold{y}^{\bold{d}_{\Bar{b}d,c\Bar{c}}}\bold{y}^{\bold{d}_{ac,\Bar{b}d}}F_{(a,\Bar{b})}=-\bold{y}^{\bold{d}_{a\ast,\Bar{b}\ast}}F_{(a,\Bar{b})}.
\end{equation}

So, we obtain 
\begin{equation*}
   F_{\mathcal{G}_{ab}}=F_{\gamma_1}F_{\gamma_2}-\bold{y}^{\bold{d}_{\gamma_1,\gamma_2}}F_{(a,\Bar{b})}=F_{ab}^B. 
\end{equation*}
\item [2)] The case in which $\tau_{[n-1]}$ is in the minimal perfect matching of $\hat{\mathcal{G}}_{\gamma_1}$ is analogous, exchanging the roles of $\gamma_1$ and $\gamma_2$. 

\end{itemize}

    \end{itemize}

\end{proof}
\begin{lemma}\label{lemma2_B}
   Let $T=\{ \tau_1, \dots, \tau_n=d, \dots, \tau_{2n-1} \}$ be a $\theta$-invariant triangulation of $\mathbf{P}_{2n+2}$ such that $\tau_n$ and $\tau_{n-1}$ are edges of a triangle of $T$ whose third edge is a boundary edge, and $\tau_n=d$ is oriented. For any $\theta$-orbit $[a,b]$ of $\mathbf{P}_{2n+2}$, $\bold{g}_{\mathcal{G}_{ab}}=g^B_{ab}$ (cf. Definition \ref{def_type_B}). 
\end{lemma}

\begin{proof}
 If $\text{Res}([a,b])=\{ \gamma \}$, by construction, an edge with label $n$ is in $\mathcal{G}_\gamma$ if and only if two edges with label $n$ are in $\hat{\mathcal{G}}_\gamma$. Therefore, 
 
\begin{itemize}
    \item If $\gamma$ does not cross $\tau_n$, $\bold{g}_{\mathcal{G}_{ab}}=\bold{g}_{\hat{\mathcal{G}}_\gamma}=D\bold{g}_{\mathcal{G}_\gamma}$.
    \item Otherwise, if $\gamma$ crosses $\tau_n$, then $\bold{g}_{\mathcal{G}_{ab}}=\bold{g}_{\hat{\mathcal{G}}_\gamma}=D\bold{g}_{\mathcal{G}_\gamma}+\bold{e}_n$, since in $D\bold{g}_{\mathcal{G}_\gamma}$ we have subtracted $\bold{e}_n$ twice, so we have to add it once.
\end{itemize} 
If $\text{Res}([a,b])=\{ \gamma_1, \gamma_2 \}$, the statement follows since the minimal matching of $\mathcal{G}_{ab}$ is the gluing of the minimal matchings of $\hat{\mathcal{G}}_{\gamma_1}$ and $\hat{\mathcal{G}}_{\gamma_2}$.    
\end{proof}

\begin{theorem}\label{theorem 3}
Let $T=\{ \tau_1, \dots, \tau_n=d, \dots, \tau_{2n-1} \}$ be a $\theta$-invariant triangulation of $\mathbf{P}_{2n+2}$ with oriented diameter $\tau_n=d$, such that $\tau_n$ and $\tau_{n-1}$ are edges of a triangle of $T$ whose third edge is a boundary edge. Let $\mathcal{A}=\mathcal{A}(T)^B$ be the cluster algebra of type $B_n$ with principal coefficients in $T$. 
    Let $[a,b]$ be an orbit of the action of $\theta$ on the diagonals of the polygon, and $x_{ab}$ the cluster variable of $\mathcal{A}$ which corresponds to $[a,b]$. Let $F_{ab}$ and $\bold{g}_{ab}$ be the $F$-polynomial and the $\bold{g}$-vector of $x_{ab}$, respectively. 
    Then $F_{ab}=F_{\mathcal{G}_{ab}}$ and $\bold{g}_{ab}=\bold{g}_{\mathcal{G}_{ab}}$.
\end{theorem}
\begin{proof}
  The result follows directly from Theorem \ref{theorem1}, Lemma \ref{lemma1_B} and Lemma \ref{lemma2_B}.  
\end{proof}
\begin{remark}
    Theorem \ref{theorem 3} extends the result of \cite{M} for cluster algebras of type $B$ to every seed whose cluster corresponds to a $\theta$-invariant triangulation $T=\{ \tau_1, \dots, \tau_n=d, \dots, \tau_{2n-1} \}$ of $\mathbf{P}_{2n+2}$, such that $\tau_n=d$ and $\tau_{n-1}$ are edges of a triangle of $T$ whose third edge is a boundary edge. 
\end{remark}
\begin{example}
    Let $[a,b]$ be the $\theta$-orbit in the triangulated octagon in Figure \ref{ex_lab_snake_graph2_}. It follows from Theorem \ref{theorem 3} that the Laurent polynomial $F_{\mathcal{G}_{ab}}$, and the integer vector $\bold{g}_{\mathcal{G}_{ab}}$, computed in Example \ref{ex_lab_snake_graph2_B}, are the $F$-polynomial, and the $\bold{g}$-vector respectively, of $x_{ab}\in \mathcal{A}(T)^B$, where $T$ is the $\theta$-invariant triangulation of the octagon in Figure \ref{ex_lab_snake_graph2_}.
\end{example}
\subsection{Type C}

\begin{definition}
    Let $\Bar{T}=\{ \tau_1, \dots, \tau_n \}$ be a triangulation of $\mathbf{P}_{n+3}$, such that $\tau_n$ is an edge of a triangle of $\Bar{T}$ whose other two edges are boundary edges. Let $\gamma$ be a diagonal of $\mathbf{P}_{n+3}$ which is not in $\Bar{T}$. We define the $labeled$ $modified$ $snake$ $graph$ $\hat{\mathcal{G}}_{\gamma}$ associated with $\gamma$ as the usual labeled snake graph $\mathcal{G}_{\gamma}$ of Definition \ref{def_snake_graph} with the following additional labels on the tile $\hat{G}_n$ with label $n$: if $l$ is a label of an edge $e$ of $G_n$, $l$ is also a label of the edge of $\hat{G}_n$ opposite to $e$.     
\end{definition}
\begin{remark}
    A cluster algebra of type $C_n$ can also be realized as a disk with one orbifold point of weight $\frac{1}{2}$, and $n+1$ boundary marked points \cite{FeST}. In \cite{canakciTumarkin}, \text{\c{C}anak\c{c}\i} and Tumarkin introduce snake and band graphs associated to curves in a triangulated orbifold with orbifold points of weight $\frac{1}{2}$, including type $C$. The tile $\hat{G}_n$ of $\hat{\mathcal{G}}_\gamma$ is the same as the tile they associate to the pending arc, i.e. the arc of the triangulation of the orbifold connecting a boundary point to the orbifold point. 
\end{remark}
\begin{example}
    In the example for $n=3$ in Figure \ref{fig:snake_graph}, we compute snake graphs $\mathcal{G}_\gamma$ and $\hat{\mathcal{G}}_\gamma$ of a diagonal $\gamma$ in a triangulated hexagon.
\begin{figure}[H]
    \centering
   \begin{tikzpicture}[scale=0.9]
    \draw (0:3cm) -- (60:3cm) -- (120:3cm) -- (180:3cm) -- (240:3cm) -- (300:3cm) -- cycle;
                    \draw (0:3cm) -- node[midway, above, xshift=-1mm, yshift=-1mm] {3} (120:3cm);
                    \draw (0:3cm) -- node[midway, above left,xshift=1.2mm] {2} (180:3cm);
                    \draw (300:3cm) -- node[midway, above left,xshift=1mm] {1} (180:3cm);
                    \draw[blue, line width=0.3mm] (60:3cm) -- node[left, yshift=-9mm] {$\gamma$} (240:3cm);
   \begin{scope}[xshift=7cm]

   \draw (-0.5,0) -- (-0.5,1) -- (0.5,1) -- (0.5,0) -- cycle;

\begin{scope}
    \draw (-0.5,0) --node[midway] {$[1]$} (-0.5,1) --node[midway, above left,xshift=2mm] {3} (0.5,1) -- node[midway] {$[2]$}(0.5,0) -- node[midway, below,xshift=1mm] {1} cycle;
     \node at (0,0.5) {2};

     \node at (-2.5,0.5) {$\mathcal{G}_\gamma=$};
     
\end{scope}

\begin{scope}[xshift=-1cm]
    \draw (-0.5,0) -- (-0.5,1) -- node[midway, above] {2}(0.5,1) -- (0.5,0) --  cycle;
    \node at (0,0.5) {1};
    
\end{scope}  

\begin{scope}[xshift=1cm]
    \draw (-0.5,0) -- (-0.5,1) --(0.5,1) -- (0.5,0) -- node[midway, below] {2}cycle;
    \node at (0,0.5) {3};
    
\end{scope} 
  \end{scope}

     \begin{scope}[xshift=7cm,yshift=-2.5cm]

   \draw (-0.5,0) -- (-0.5,1) -- (0.5,1) -- (0.5,0) -- cycle;

\begin{scope}
    \draw (-0.5,0) --node[midway] {$[1]$} (-0.5,1) --node[midway, above left,xshift=2mm] {3} (0.5,1) -- node[midway] {$[2]$}(0.5,0) -- node[midway, below,xshift=1mm] {1} cycle;
     \node at (0,0.5) {2};

     \node at (-2.5,0.5) {$\hat{\mathcal{G}_\gamma}=$};
     
\end{scope}

\begin{scope}[xshift=-1cm]
    \draw (-0.5,0) -- (-0.5,1) -- node[midway, above] {2}(0.5,1) -- (0.5,0) --  cycle;
    \node at (0,0.5) {1};
    
\end{scope}  

\begin{scope}[xshift=1cm]
    \draw (-0.5,0) -- (-0.5,1) --node[midway, above] {2}(0.5,1) -- node[midway] {$[2]$}(0.5,0) -- node[midway, below] {2}cycle;
    \node at (0,0.5) {3};
    
\end{scope} 
  \end{scope}
   \end{tikzpicture}
    \caption{The snake graphs $\mathcal{G}_\gamma$ and $\hat{\mathcal{G}_\gamma}$ for a diagonal $\gamma$ in a triangulated hexagon (type $C$).}
    \label{fig:snake_graph}
\end{figure}

\end{example}

\begin{definition}
    Let $T=\{ \tau_1, \dots, \tau_n=d, \dots, \tau_{2n-1} \}$ be a $\theta$-invariant triangulation of $\mathbf{P}_{2n+2}$ with oriented diameter $d$. Let $[a,b]$ be a $\theta$-orbit which is not in $T$. We associate to $[a,b]$ the labeled modified snake graph $\mathcal{G}_{ab}$ defined in the following way:
    \begin{itemize}
        \item if $\Tilde{\text{Res}}([a,b])=\{ \Tilde{\gamma} \}$, then $\mathcal{G}_{ab}:=\hat{\mathcal{G}}_{\tilde\gamma}$;
        \item if $\Tilde{\text{Res}}([a,b])=\{ \Tilde{\gamma_1},\Tilde{\gamma_2} \}$, then $\mathcal{G}_{ab}$ is obtained by gluing $\hat{\mathcal{G}}_{\tilde\gamma_1}$ and $\hat{\mathcal{G}}_{\tilde\gamma_2}$ along their common exterior edge.
    \end{itemize}
\end{definition}

\begin{remark}\label{rmk_gluing_edge}
    In the case where $\Tilde{\text{Res}}([a,b])=\{ \Tilde{\gamma_1},\Tilde{\gamma_2} \}$, with the notation of Definition \ref{def_type_C}, the edge along which we glue $\hat{\mathcal{G}}_{\tilde\gamma_1}$ and $\hat{\mathcal{G}}_{\tilde\gamma_2}$ is $(\Bar{b},c)$ if $[a,b]=[a,\Bar{a}]$ is a diameter, while it is $(\Bar{c},\Bar{d})$ if $[a,b]$ is a pair of diagonals which cross $d$.
\end{remark}

\begin{remark}
    Let $[a,b]$ be a $\theta$-orbit such that $\text{Res}([a,b])=\{ \gamma_1,\gamma_2 \}$. Then $\mathcal{G}_{ab}$ is obtained by superimposing $\mathcal{G}_{\gamma_1}$ and $\mathcal{G}_{\gamma_2}$ over their tile $G_n$ with label $n$, in the only way such that $G_n$ has different relative orientation with respect to $\Bar{T}=\text{Res}(T)$ in $\mathcal{G}_{\gamma_1}$ and $\mathcal{G}_{\gamma_2}$.
\end{remark}

We define the perfect matching polynomial $F_{\mathcal{G}_{ab}}$ and the $\bold{g}$-vector $\bold{g}_{\mathcal{G}_{ab}}$ of $\mathcal{G}_{ab}$ as in Definition \ref{def_f_poly_lab_mod_sg}, where for the height monomial we can use Definition \ref{def_h(P)_original} since, unlike type $B$, we do not have an additional edge. 

\begin{example}\label{ex_lab_snake_graph_C}
We compute the labeled modified snake graph $\mathcal{G}_{ab}$ of the $\theta$-orbit $[a,b]$ in Figure \ref{ex_lab_snake_graph}.
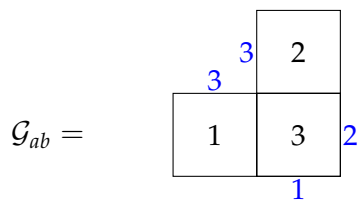
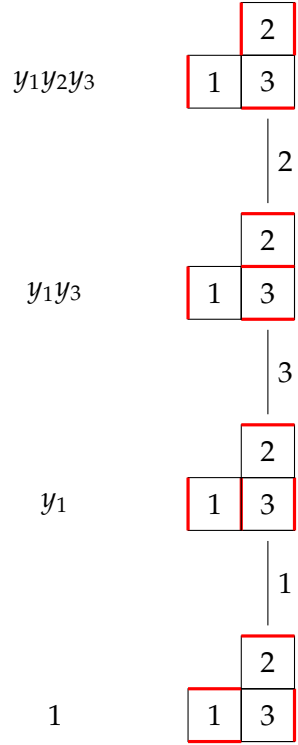
\begin{figure}[H]
    \centering
   \begin{tikzpicture}[scale=1.1]
  
 
   \draw (-0.5,0) -- (-0.5,1) -- (0.5,1) -- (0.5,0) -- cycle;

\node at (-3,0.5) {$\mathcal{G}_{ab}=$};

\begin{scope}
    \filldraw[fill=red!10] (-0.5,0) --node[midway] {} (-0.5,1) --node[midway, above left,xshift=2mm] {} (0.5,1) -- node[midway,right,xshift=-1mm] {\textcolor{blue}{2}}(0.5,0) -- node[midway, below,yshift=1mm] {\textcolor{blue}{1}} cycle;
     \node at (0,0.5) {3};

\end{scope}

\begin{scope}[xshift=-1cm]
    \filldraw[fill=blue!10] (-0.5,0) -- (-0.5,1) -- node[midway, above,yshift=-1mm] {\textcolor{blue}{3}}(0.5,1) -- (0.5,0) --  cycle;
    \node at (0,0.5) {1};
\end{scope}  

\begin{scope}[yshift=1cm]
    \filldraw[fill=red!10] (-0.5,0) --node[midway, left,xshift=1mm] {\textcolor{blue}{3}} (-0.5,1) --(0.5,1) -- (0.5,0) -- cycle;
    \node at (0,0.5) {2};
\end{scope} 
\end{tikzpicture}
\end{figure}
\begin{figure}
    \centering
     \begin{tikzpicture}[scale=3.1]
  
 
   \draw (-0.5,0) -- (-0.5,1) -- (0.5,1) -- (0.5,0) -- cycle;
\node at (-0.5,0) {$\bullet$};
\node at (-0.5,1) {$\bullet$};
\node at (0.5,0) {$\bullet$};
\node at (0.5,1) {$\bullet$};
\node at (0,0.5) {$\bold\times$};
\draw (0,0.5) -- node[midway]{3} (0.5,1);
\draw (-0.5,0) arc  (180:90:1)node[midway]{1};
\draw (0.5,1) arc (00:-90:1)node[midway]{2};
\draw[magenta] (0.5,0) arc (00:90:1)node[midway, left,xshift=-3mm,yshift=2mm]{$\gamma$};

\end{tikzpicture}
    \caption{A triangulated orbifold with one orbifold point of weight $\frac{1}{2}$.}
    \label{triang_orbif}
\end{figure}
\begin{remark}
    $\mathcal{G}_{ab}$ is the snake graph associated in \cite{canakciTumarkin} to the arc $\gamma$ in the triangulated orbifold in Figure \ref{triang_orbif}.
\end{remark}

\begin{figure}
                               \centering
                \begin{tikzpicture}[scale=0.5]
                    \draw (90:3cm) -- (135:3cm) -- (180:3cm) -- (225:3cm) -- (270:3cm) -- (315:3cm) -- (360:3cm) -- (45:3cm) -- cycle;
                    \draw (270:3cm) -- (180:3cm);;
                    \draw (90:3cm) -- (360:3cm);
                    \draw[-{Latex[length=2mm]}] (270:3cm) -- node[midway, above left,xshift=1mm] {} (90:3cm);
                    \draw (180:3cm) -- node[midway, above left,xshift=1mm] {} (90:3cm);
                    \draw (360:3cm) -- node[midway, above left,xshift=1mm] {} (270:3cm);
                    \draw[cyan, line width=0.3mm] (135:3cm) -- (45:3cm); 
                    \draw[cyan, line width=0.3mm] (225:3cm) -- (315:3cm);

                    \node at (135:3cm) [left] {$\Bar{b}$};
                    \node at (225:3cm) [left] {$a$};
                    \node at (45:3cm) [right] {$\Bar{a}$};
                    \node at (315:3cm) [right] {$b$};

                                \begin{scope}[xshift=8cm]
                      \draw (90:3cm) -- (135:3cm) -- (180:3cm) -- (225:3cm) -- (270:3cm) -- (360:3cm) -- cycle;
 
                    \draw (180:3cm) -- node[midway, above left,xshift=1mm] {2}(270:3cm);;
                    \draw (90:3cm) -- node[midway, above left,xshift=1mm] {3} (270:3cm);
                    \draw (180:3cm) -- node[midway, left,xshift=1mm] {1} (90:3cm);

\draw[red, line width=0.3mm] (225:3cm) -- node[midway,below] {$\tilde\gamma_1$} (360:3cm); 
\draw[blue, line width=0.3mm] (135:3cm) -- node[midway, right,xshift=-1mm] {$\tilde\gamma_2$} (270:3cm);

\node at (135:3cm) [left] {$\Bar{b}$};
\node at (90:3cm) [above] {$c$};
\node at (270:3cm) [below] {$\Bar{c}$};
\node at (180:3cm) [left] {$\Bar{d}$};
\node at (225:3cm) [left] {$a$};
\node at (360:3cm) [right] {$\ast$}; 
                    \end{scope}
                \end{tikzpicture}
                               \caption{A $\theta$-orbit $[a,b]$ in a triangulated octagon and its rotated restriction (type $C_3$).}
                               \label{ex_lab_snake_graph}
                           \end{figure}

\begin{figure}[H]
    \centering
   \begin{tikzpicture}[scale=0.7]
  
 
   \draw (-0.5,0) -- (-0.5,1) -- (0.5,1) -- (0.5,0) -- cycle;

\node at (-4,0.5) {$y_1y_2y_3$};

\begin{scope}
    \draw (-0.5,0) --node[midway] {} (-0.5,1) --node[midway, above left,xshift=2mm] {} (0.5,1) -- node[midway] {}(0.5,0) -- node[midway, below,xshift=1mm] {} cycle;
    \draw[red,very thick] (-0.5,0) -- (0.5,0);
     \node at (0,0.5) {3};

\end{scope}

\begin{scope}[xshift=-1cm]
    \draw (-0.5,0) -- (-0.5,1) -- node[midway, above] {}(0.5,1) -- (0.5,0) --  cycle;
    \node at (0,0.5) {1};
   \draw[red,very thick] (-0.5,0) -- (-0.5,1); 
\end{scope}  

\begin{scope}[yshift=1cm]
    \draw (-0.5,0) -- (-0.5,1) --(0.5,1) -- (0.5,0) -- node[midway, below] {}cycle;
    \node at (0,0.5) {2};
   \draw[red,very thick] (-0.5,0) -- (-0.5,1);
   \draw[red,very thick] (0.5,0) -- (0.5,1);
\end{scope} 
\draw (0,-0.2) -- node[midway, right] {2} (0,-1.8);
\begin{scope}[yshift=-4cm]
      \draw (-0.5,0) -- (-0.5,1) -- (0.5,1) -- (0.5,0) -- cycle;
\draw (0,-0.2) -- node[midway,right] {3}  (0,-1.8);
\node at (-4,0.5) {$y_1y_3$};

\begin{scope}
    \draw (-0.5,0) --node[midway] {} (-0.5,1) --node[midway, above left,xshift=2mm] {} (0.5,1) -- node[midway] {}(0.5,0) -- node[midway, below,xshift=1mm] {} cycle;
   \draw[red,very thick] (-0.5,0) -- (0.5,0);
     \node at (0,0.5) {3};

\end{scope}

\begin{scope}[xshift=-1cm]
    \draw (-0.5,0) -- (-0.5,1) -- node[midway, above] {}(0.5,1) -- (0.5,0) --  cycle;
    \node at (0,0.5) {1};
   \draw[red,very thick] (-0.5,0) -- (-0.5,1); 
\end{scope}  

\begin{scope}[yshift=1cm]
    \draw (-0.5,0) -- (-0.5,1) --(0.5,1) -- (0.5,0) -- node[midway, below] {}cycle;
    \node at (0,0.5) {2};
    \draw[red,very thick] (-0.5,0) -- (0.5,0);
    \draw[red,very thick] (-0.5,1) -- (0.5,1);
\end{scope}
\end{scope}
\begin{scope}[yshift=-8cm]
      \draw (-0.5,0) -- (-0.5,1) -- (0.5,1) -- (0.5,0) -- cycle;
\node at (-4,0.5) {$y_1$};
\draw (0,-0.2) -- node[midway,right] {1} (0,-1.8);

\begin{scope}
    \draw (-0.5,0) --node[midway] {} (-0.5,1) --node[midway, above left,xshift=2mm] {} (0.5,1) -- node[midway] {}(0.5,0) -- node[midway, below,xshift=1mm] {} cycle;
     \node at (0,0.5) {3};
\draw[red,very thick] (-0.5,0) -- (-0.5,1);
 \draw[red,very thick] (0.5,0) -- (0.5,1);  
     
\end{scope}

\begin{scope}[xshift=-1cm]
    \draw (-0.5,0) -- (-0.5,1) -- node[midway, above] {}(0.5,1) -- (0.5,0) --  cycle;
    \node at (0,0.5) {1};
\draw[red,very thick] (-0.5,0) -- (-0.5,1);    
\end{scope}  

\begin{scope}[yshift=1cm]
    \draw (-0.5,0) -- (-0.5,1) --(0.5,1) -- (0.5,0) -- node[midway, below] {}cycle;
    \node at (0,0.5) {2};
 \draw[red,very thick] (-0.5,1) -- (0.5,1);   
\end{scope}
\end{scope}
\begin{scope}[yshift=-12cm]
      \draw (-0.5,0) -- (-0.5,1) -- (0.5,1) -- (0.5,0) -- cycle;
\node at (-4,0.5) {$1$};

\begin{scope}
    \draw (-0.5,0) --node[midway] {} (-0.5,1) --node[midway, above left,xshift=2mm] {} (0.5,1) -- node[midway] {}(0.5,0) -- node[midway, below,xshift=1mm] {} cycle;
     \node at (0,0.5) {3};
\draw[red,very thick] (0.5,0) -- (0.5,1);

\end{scope}

\begin{scope}[xshift=-1cm]
    \draw (-0.5,0) -- (-0.5,1) -- node[midway, above] {}(0.5,1) -- (0.5,0) --  cycle;
    \node at (0,0.5) {1};
 \draw[red,very thick] (-0.5,0) -- (0.5,0);
 \draw[red,very thick] (-0.5,1) -- (0.5,1);
\end{scope}  

\begin{scope}[yshift=1cm]
    \draw (-0.5,0) -- (-0.5,1) --(0.5,1) -- (0.5,0) -- node[midway, below] {}cycle;
    \node at (0,0.5) {2};
\draw[red,very thick] (-0.5,1) -- (0.5,1);    
\end{scope}
\end{scope}    
   \end{tikzpicture}
    \caption{The poset of perfect matchings of $\mathcal{G}_{ab}$, and the corresponding monomials.}
    \label{set_pm_C}
\end{figure} 
Moreover, we compute the perfect matching polynomial $F_{\mathcal{G}_{ab}}$ and the $\bold{g}$-vector $\bold{g}_{\mathcal{G}_{ab}}$ of $\mathcal{G}_{ab}$. The set of all perfect matchings of $\mathcal{G}_{ab}$, with the corresponding monomials, is illustrated in Figure \ref{set_pm_C}. We have 
\begin{align*}
    F_{\mathcal{G}_{ab}}=y_1y_2y_3+y_1y_3+y_1+1,
\end{align*}
and
\begin{align*}
    \bold{g}_{\mathcal{G}_{ab}}=\begin{pmatrix} 0\\ 1\\ 1
    \end{pmatrix}-\begin{pmatrix} 1\\ 1\\1
    \end{pmatrix}=\begin{pmatrix} -1\\ 0\\ 0
    \end{pmatrix}.
\end{align*}
\end{example}

\begin{lemma}\label{lemma1_C}
   Let $T=\{ \tau_1, \dots, \tau_n=d, \dots, \tau_{2n-1} \}$ be a $\theta$-invariant triangulation of $\mathbf{P}_{2n+2}$ with oriented diameter $d$. For any $\theta$-orbit $[a,b]$ of $\mathbf{P}_{2n+2}$, $F_{\mathcal{G}_{ab}}=F^C_{ab}$ (cf. Definition \ref{def_type_C}).  
\end{lemma}
\begin{proof}
 If $\Tilde{\text{Res}}([a,b])=\{ \Tilde{\gamma} \}$, the statement holds since $F_{\mathcal{G}_{ab}}=F_{\hat{\mathcal{G}}_{\tilde\gamma}}=F_{\mathcal{G}_{\tilde\gamma}}$. Otherwise, if $\Tilde{\text{Res}}([a,b])=\{ \Tilde{\gamma_1},\Tilde{\gamma_2} \}$, then
\begin{equation}
    F_{\mathcal{G}_{ab}}=F_{\hat{\mathcal{G}}_{\tilde\gamma_1}}F_{\hat{\mathcal{G}}_{\tilde\gamma_2}}-R=F_{\mathcal{G}_{\tilde\gamma_1}}F_{\mathcal{G}_{\tilde\gamma_2}}-R,
\end{equation}
where $R$ is the sum of the monomials which correspond to gluing of perfect matchings of $\hat{\mathcal{G}}_{\tilde\gamma_1}$ (in red in Figure \ref{forbidden}) and perfect matchings of $\hat{\mathcal{G}}_{\tilde\gamma_2}$ (in blue in Figure \ref{forbidden}) which are not perfect matchings of $\mathcal{G}_{ab}$. They are all of the form
\begin{figure}[H]
    \centering
   \begin{tikzpicture}[scale=1.6]
   
    \filldraw[fill=red!10] (-0.5,0) -- node[midway, left] {$\cdots$} (-0.5,1) -- node[midway, above] {}(0.5,1) --node[midway,yshift=0mm] {} (0.5,0) -- node[midway, below,yshift=0mm] {$\vdots$} cycle;
    \node at (0,0.5) {$n$};
    \draw[red, very thick] (-0.5,1) -- node[above,yshift=-1mm] {\textcolor{black}{$e_1$}}(0.5,1);
    \draw[red, very thick] (0.5,0) --node[below,yshift=1mm] {\textcolor{black}{$e_2$}}(-0.5,0);
\begin{scope}[xshift=1cm]
    \filldraw[fill=blue!10] (-0.5,0) -- (-0.5,1) --node[midway, above,yshift=1mm] {$\vdots$}(0.5,1) -- node[midway,right] {$\cdots$}(0.5,0) -- node[midway, below] {}cycle;
    \node at (0,0.5) {$i$};
    \draw[red, very thick] (-0.5,1) --node[above,yshift=-1mm] {\textcolor{black}{$e_3$}}(0.5,1);
    \draw[red, very thick] (0.5,0) --node[below,yshift=1mm] {\textcolor{black}{$e_4$}}(-0.5,0);
    
  \end{scope}  

\draw (0.5,1) --node[midway,xshift=0.5mm] {$\tau_{[i]}$} (0.5,0);
   \end{tikzpicture}
                               \label{forbidden}
   
\end{figure}

Therefore, we have to describe all perfect matchings of $\hat{\mathcal{G}}_{\tilde\gamma_1}$ and of $\hat{\mathcal{G}}_{\tilde\gamma_2}$ which do not contain the edge with label $\tau_{[i]}$, along which we glue $\hat{G}_n$ of $\hat{\mathcal{G}}_{\tilde\gamma_1}$ and $\hat{G}_i$ of $\hat{\mathcal{G}}_{\tilde\gamma_2}$. 

We prove the statement in the case where $[a,b]=[a,\Bar{a}]$ is a diameter. If $[a,b]$ is a pair of diagonals which cross $d$, the proof is completely analogous. We consider the type $A$ exchange relation corresponding to the crossing of diagonals $\tilde\gamma_1$ and $\tau_n$.

\begin{figure}[H]
                               \centering
                \begin{tikzpicture}[scale=0.5]
                     \draw (90:3cm) -- (110:3cm) -- (130:3cm) -- (150:3cm)-- (170:3cm) -- (190:3cm) -- (210:3cm)--(230:3cm)--(250:3cm) -- (270:3cm) -- (360:3cm)  --  cycle;
                    \draw (90:3cm) -- node[midway, above left,xshift=1mm] {} (270:3cm);
                   
                    \draw[cyan, line width=0.3mm] (90:3cm) -- node[midway, xshift=-2mm] {$\tau_n$}(270:3cm);
                     \draw[red, line width=0.3mm] (150:3cm) -- node[midway, above left, xshift=-1mm] {$\tilde\gamma_1$}(360:3cm);

                    \node at (90:3cm) [above] {$b$};
                    \node at (270:3cm) [below] {$\Bar{b}$};

\node at (150:3cm) [left] {$a$};
                    \node at (360:3cm) [right] {$\ast$};
                      
                \end{tikzpicture}
                               \caption{Type $A$ exchange relation corresponding to the crossing of $\tilde\gamma_1$ and $\tau_n$.}
                               \label{fig:enter-label}
                           \end{figure}  

We have that
\begin{equation}\label{eq_1}
    F_{\tilde\gamma_1}=\bold{y}^{\bold{d}_{a\Bar{b},b\ast}}F_{(a,b)}+F_{(a,\Bar{b})}.
\end{equation}
At this point, we have two cases to consider.
\begin{itemize}
    \item[1)] $\tau_{[i]}$ is not in the minimal perfect matching of $\hat{\mathcal{G}}_{\tilde\gamma_1}$. So the red edges $e_1$, $e_2$ of $\hat{G}_n$ in Figure \ref{forbidden}
  are in the minimal perfect matching of $\hat{\mathcal{G}}_{\tilde\gamma_1}$. So the sum of the monomials which correspond to the perfect matchings of $\hat{\mathcal{G}}_{\tilde\gamma_1}$ which contain them in the right hand side of \ref{eq_1} is $F_{(a,\Bar{b})}$. Moreover, the fact that $e_1$, $e_2$ are in the minimal matching of $\hat{\mathcal{G}}_{\tilde\gamma_1}$ means that the monomial $y_n$ is a summand of $F_{\mathcal{G}_{\tilde\gamma_1}}$. Therefore, $\tau_i$ has to be clockwise from $\tau_n$.

We consider the type $A$ exchange relation corresponding to the crossing of diagonals $\tilde\gamma_2$ and $\tau_i$.

\begin{figure}[H]
                               \centering
                \begin{tikzpicture}[scale=0.5]
                       \draw (90:3cm) -- (110:3cm) -- (130:3cm) -- (150:3cm)-- (170:3cm) -- (190:3cm) -- (210:3cm)--(230:3cm)--(250:3cm) -- (270:3cm) -- (360:3cm)  --  cycle;
                      \draw (90:3cm) -- node[midway, right,xshift=-1mm] {$\tau_n$} (270:3cm);
                   
                    \draw[cyan, line width=0.3mm] (110:3cm) -- node[midway, xshift=-2mm] {$\tau_i$}(270:3cm);
                     \draw[blue, line width=0.3mm] (150:3cm) -- node[midway, below right] {$\tilde\gamma_2$}(90:3cm);

                    \node at (90:3cm) [above] {$b$};
                    \node at (270:3cm) [below] {$\Bar{b}$};

\node at (150:3cm) [left] {$a$};
\node at (120:3cm) [above] {$c$};
                    \node at (360:3cm) [right] {$\ast$};
                      
                \end{tikzpicture}
                               \caption{Type $A$ exchange relation corresponding to the crossing of $\tilde\gamma_2$ and $\tau_i$.}
                               \label{fig:enter-label}
                           \end{figure} 
We have that
\begin{equation}\label{eq_22}
    F_{\tilde\gamma_2}=\bold{y}^{\bold{d}_{bc,a\Bar{b}}}F_{(a,c)}+F_{(a,\Bar{b})}.
\end{equation}
Since $e_1$, $e_2$ are in the minimal perfect matching of $\hat{\mathcal{G}}_{\tilde\gamma_1}$, the red edges $e_3$, $e_4$ of $\hat{G}_i$ in Figure \ref{forbidden} cannot be in the minimal perfect matching of $\hat{\mathcal{G}}_{\tilde\gamma_2}$. So the sum of the monomials which correspond to the perfect matchings of $\hat{\mathcal{G}}_{\tilde\gamma_2}$ which contain $e_3$, $e_4$ in the right hand side of \ref{eq_22} is $\bold{y}^{\bold{d}_{bc,a\Bar{b}}}F_{(a,c)}$. 
Therefore,
\begin{equation}
    R=\bold{y}^{\bold{d}_{bc,a\Bar{b}}}F_{(a,\Bar{b})}F_{(a,c)}=\bold{y}^{\bold{d}_{bc,a\ast}}F_{(a,\Bar{b})}F_{(a,c)}.
\end{equation}

So, we obtain (cf. Remark \ref{rmk_coeff_typeC})
\begin{equation*}
    F_{\mathcal{G}_{ab}}=F_{\tilde\gamma_1}F_{\tilde\gamma_2}-\bold{y}^{\bold{d}_{bc,a\ast}}F_{(a,\Bar{b})}F_{(a,c)}=F_{ab}^C.
\end{equation*}

\item [2)] $\tau_{[i]}$ is in the minimal perfect matching of $\hat{\mathcal{G}}_{\tilde\gamma_1}$. So $e_1$, $e_2$
  are not in the minimal perfect matching of $\hat{\mathcal{G}}_{\tilde\gamma_1}$. So the sum of the monomials which correspond to the perfect matchings of $\hat{\mathcal{G}}_{\tilde\gamma_1}$ which contain $e_1$, $e_2$ in the right hand side of \ref{eq_1} is $\bold{y}^{\bold{d}_{a\Bar{b},b\ast}}F_{(a,b)}$. Moreover, the fact that $e_1$, $e_2$ are not in the minimal matching of $\hat{\mathcal{G}}_{\tilde\gamma_1}$ means that the monomial $y_n$ is not a summand of $F_{\mathcal{G}_{\tilde\gamma_1}}$. Therefore, $\tau_i$ has to be counterclockwise from $\tau_n$.

We consider the type $A$ exchange relation corresponding to the crossing of diagonals $\tilde\gamma_2$ and $\tau_i$.

\begin{figure}[H]
                               \centering
                \begin{tikzpicture}[scale=0.5]
                       \draw (90:3cm) -- (110:3cm) -- (130:3cm) -- (150:3cm)-- (170:3cm) -- (190:3cm) -- (210:3cm)--(230:3cm)--(250:3cm) -- (270:3cm) -- (360:3cm)  --  cycle;
                     \draw (90:3cm) -- node[midway, right,xshift=-1mm] {$\tau_n$} (270:3cm);
                   
                    \draw[cyan, line width=0.3mm] (90:3cm) -- node[midway, xshift=-2mm] {$\tau_i$}(250:3cm);
                     \draw[blue, line width=0.3mm] (150:3cm) -- node[midway, right,yshift=-2mm] {$\tilde\gamma_2$}(270:3cm);

                    \node at (90:3cm) [above] {$b$};
                    \node at (270:3cm) [below] {$\Bar{b}$};

\node at (150:3cm) [left] {$a$};
\node at (240:3cm) [below] {$c$};
                    \node at (360:3cm) [right] {$\ast$};
                      
                \end{tikzpicture}
                               \caption{Type $A$ exchange relation corresponding to the crossing of $\tilde\gamma_2$ and $\tau_i$.}
                               \label{fig:enter-label}
                           \end{figure} 
We have that
\begin{equation}\label{eq_2}
    F_{\tilde\gamma_2}=\bold{y}^{\bold{d}_{ac,b\Bar{b}}}F_{(a,b)}+F_{(a,c)}.
\end{equation}
Since $e_1$, $e_2$ are not in the minimal perfect matching of $\hat{\mathcal{G}}_{\tilde\gamma_1}$, $e_3$, $e_4$ must be in the minimal perfect matching of $\hat{\mathcal{G}}_{\tilde\gamma_2}$. So the sum of the monomials which correspond to the perfect matchings of $\hat{\mathcal{G}}_{\tilde\gamma_2}$ which contain $e_3$, $e_4$ in the right hand side of \ref{eq_2} is $F_{(a,c)}$. 
Therefore,
\begin{equation}
    R=\bold{y}^{\bold{d}_{a\Bar{b},b\ast}}F_{(a,b)}F_{(a,c)}.
\end{equation}

So, we obtain (cf. Remark \ref{rmk_coeff_typeC})
\begin{equation*}
    F_{\mathcal{G}_{ab}}=F_{\tilde\gamma_1}F_{\tilde\gamma_2}-\bold{y}^{\bold{d}_{a\Bar{b},b\ast}}F_{(a,b)}F_{(a,c)}=F_{ab}^C.
\end{equation*}
\end{itemize}
\end{proof}

\begin{lemma}\label{lemma2_C}
    Let $T=\{ \tau_1, \dots, \tau_n=d, \dots, \tau_{2n-1} \}$ be a $\theta$-invariant triangulation of $\mathbf{P}_{2n+2}$ with oriented diameter $d$. For any $\theta$-orbit $[a,b]$ of $\mathbf{P}_{2n+2}$, $\bold{g}_{\mathcal{G}_{ab}}=g^C_{ab}$ (cf. Definition \ref{def_type_C}).
\end{lemma}

\begin{proof}
If $\Tilde{\text{Res}}([a,b])=\{\tilde\gamma \}$, we have two cases to consider.
\begin{itemize}
    \item If $\tau_i$ and $\tau_n$ are two different sides of a triangle of $T$, $\tau_i$ is clockwise from $\tau_n$, and $\tilde\gamma$ crosses $\tau_n$, then the edge of $\hat{G}_n$ with label $i$ and its opposite must be in the minimal perfect matching $P_-(\mathcal{G}_{ab})$ of $\mathcal{G}_{ab}$. Since they both have label $i$ in $\hat{G}_n$, it follows that $\bold{g}_{\mathcal{G}_{ab}}=\bold{g}_{\mathcal{G}_{\tilde\gamma}}+\bold{e}_i=\bold{g}_{\tilde\gamma}+\bold{e}_i$.
    \item Otherwise, $\bold{g}_{\mathcal{G}_{ab}}=\bold{g}_{\mathcal{G}_{\tilde\gamma}}=\bold{g}_{\tilde\gamma}$.
\end{itemize} 
If $\Tilde{\text{Res}}([a,b])=\{ \tilde\gamma_1, \tilde\gamma_2 \}$, the statement follows since the minimal matching of $\mathcal{G}_{ab}$ is the gluing of the minimal matchings of $\hat{\mathcal{G}}_{\tilde\gamma_1}$ and $\hat{\mathcal{G}}_{\tilde\gamma_2}$. If $\tau_i$ and $\tau_n$ are two different sides of a triangle of $T$ and $\tau_i$ is clockwise from $\tau_n$, we have to subtract the vector of the canonical basis of $\mathbb{R}^n$ which corresponds to the edge $e$ of the triangle containing $\tau_n$ along which we glue $\hat{\mathcal{G}}_{\tilde\gamma_1}$ and $\hat{\mathcal{G}}_{\tilde\gamma_2}$, i.e. $\bold{g}_{(\Bar{c},\Bar{d})}$ (resp. $\bold{g}_{(\Bar{b},c)}$ if $[a,b]=(a,\Bar{a})$ is a diameter) with the notation of Definition \ref{def_type_C} (cf. Remark \ref{rmk_gluing_edge}). We have to do this since $e$ is in the minimal perfect matching of $\hat{\mathcal{G}}_{\tilde\gamma_2}$, but it is not in $P_-(\mathcal{G}_{ab})$ since it becomes an interior edge of $\mathcal{G}_{ab}$ after gluing.
\end{proof}

\begin{theorem}\label{theorem 4}
 Let $T$ be a $\theta$-invariant triangulation of $\mathbf{P}_{2n+2}$ with oriented diameter $d$, and let $\mathcal{A}=\mathcal{A}(T)^C$ be the cluster algebra of type $C_n$ with principal coefficients in $T$. 
    Let $[a,b]$ be an orbit of the action of $\theta$ on the diagonals of the polygon, and $x_{ab}$ the cluster variable of $\mathcal{A}$ which corresponds to $[a,b]$. Let $F_{ab}$ and $\bold{g}_{ab}$ be the $F$-polynomial and the $\bold{g}$-vector of $x_{ab}$, respectively. 
    Then $F_{ab}=F_{\mathcal{G}_{ab}}$ and $\bold{g}_{ab}=\bold{g}_{\mathcal{G}_{ab}}$.
\end{theorem}
\begin{proof}
  The result follows directly from Theorem \ref{theorem 2}, Lemma \ref{lemma1_C} and Lemma \ref{lemma2_C}.  
\end{proof}
\begin{remark}
    Theorem \ref{theorem 4} extends the result of \cite{M} for cluster algebras of type $C$ to every seed. 
\end{remark}

\begin{example}
    Let $[a,b]$ be the $\theta$-orbit in the triangulated octagon in Figure \ref{ex_lab_snake_graph}. It follows from Theorem \ref{theorem 4} that the Laurent polynomial $F_{\mathcal{G}_{ab}}$, and the integer vector $\bold{g}_{\mathcal{G}_{ab}}$, computed in Example \ref{ex_lab_snake_graph_C}, are the $F$-polynomial, and the $\bold{g}$-vector respectively, of $x_{ab}\in \mathcal{A}(T)^C$, where $T$ is the $\theta$-invariant triangulation of the octagon in Figure \ref{ex_lab_snake_graph}.
\end{example}

   \begin{figure}
    \centering
    \begin{tikzpicture}[scale=0.7]

\node at (-8,0) {$y_1y_2^2y_3^2$};
\node at (-8,-2.25) {$+$};
\node at (-8,-4.5) {$y_2^2y_3^2+y_1y_2y_3^2$};
\node at (-8,-6.75) {$+$};
\node at (-8,-9) {$2y_2y_3^2$};
\node at (-8,-11.5) {$+$};
\node at (-8,-13.5) {$y_2y_3+y_3^2$};
\node at (-8,-15.75) {$+$};
\node at (-8,-18) {$2y_3$};
\node at (-8,-20.25) {$+$};
\node at (-8,-22.75) {$1$};
\node at (-8,1.75) {$+$};
\node at (-8,3.5) {$y_2y_3$};
    
    \draw (-1,0) .. controls (-6,1.5) and (0,2.3)  .. (2.5,0.85);
   \draw (0:1) -- node[midway,right] {}(60:1) -- (120:1) -- node[midway,left] {}(180:1);
    \draw (300:1) -- node[midway,below] {}(240:1);
\draw[red,very thick] (0:1) -- (60:1);
\draw[red,very thick] (240:1) -- (300:1);
\node at (0,0) {$2$};
\begin{scope}[xshift=-0.15cm, yshift=-0.5cm ,rotate=-30]
   \draw (-1,0) -- (-1,1) -- (-2,1) -- (-2,0) -- cycle;
   \draw[red,very thick] (-1,0) -- (-2,0);
   \draw[red,very thick] (-1,1) -- (-2,1);
     \node at (-1.5,0.5) {1}; 
\end{scope}
\begin{scope}[xshift=3cm]
\draw[] (0:1) -- node[pos=0.5, right, xshift=-4mm, yshift=0mm] {} (60:1);
    \draw (60:1) -- node[midway,above] {} (120:1) -- (180:1) --(240:1)-- (300:1) -- (0:1); 
    \node at (0,0) {$2$};
    \draw[red,very thick] (120:1) -- (180:1);
\draw[red,very thick] (240:1) -- (300:1);
 \draw[red,very thick] (0:1) -- (60:1);
\end{scope}

\begin{scope}[xshift=-1.5cm,yshift=-0.855cm]
    \draw[] (0:1) -- node[pos=0.5, right, xshift=-4mm, yshift=0mm] {} (60:1);
    \draw (60:1) -- node[midway,above] {} (120:1) -- (180:1) -- (0:1);
 \draw[red,very thick] (120:1) -- (180:1);
    
    \node at (0,0.4) {$3$};
    
\end{scope}
\begin{scope}[xshift=1.5cm,yshift=-0.855cm]
    \draw (0:1) -- node[midway] {}(60:1) --  (120:1); 
    \draw[] (120:1) -- node[pos=0.5, left,xshift=2mm] {}(180:1);
   \draw (180:1) -- node[midway,below] {}(0:1);
    
    \node at (0,0.4) {$3$};
    \node at (-0.3,1.6) {$I$};
    
\end{scope}
\draw (1,-1) -- node[midway,left] {1} (-1,-3);
\draw (1,-1) -- node[midway,right] {2}(3,-3);
\draw (1,2.8) -- node[midway,right] {$I$} (1,1.6);
\begin{scope}[xshift=-3cm,yshift=-4.5cm]
     \draw (-1,0) .. controls (-6,1.5) and (0,2.3)  .. (2.5,0.85);
   \draw (0:1) -- node[midway,right] {}(60:1) -- (120:1) -- node[midway,left] {}(180:1);
    \draw (300:1) -- node[midway,below] {}(240:1);
\draw[red,very thick] (0:1) -- (60:1);
\draw[red,very thick] (240:1) -- (300:1);
\node at (0,0) {$2$};
\begin{scope}[xshift=-0.15cm, yshift=-0.5cm ,rotate=-30]
   \draw (-1,0) -- (-1,1) -- (-2,1) -- (-2,0) -- cycle;
   \draw[red,very thick] (-1,0) -- (-1,1);
   \draw[red,very thick] (-2,0) -- (-2,1);
     \node at (-1.5,0.5) {1}; 
\end{scope}
\begin{scope}[xshift=3cm]
\draw[] (0:1) -- node[pos=0.5, right, xshift=-4mm, yshift=0mm] {} (60:1);
    \draw (60:1) -- node[midway,above] {} (120:1) -- (180:1) --(240:1)-- (300:1) -- (0:1); 
    \node at (0,0) {$2$};
        \draw[red,very thick] (120:1) -- (180:1);
\draw[red,very thick] (240:1) -- (300:1);
 \draw[red,very thick] (0:1) -- (60:1);
\end{scope}

\begin{scope}[xshift=-1.5cm,yshift=-0.855cm]
    \draw[] (0:1) -- node[pos=0.5, right, xshift=-4mm, yshift=0mm] {} (60:1);
    \draw (60:1) -- node[midway,above] {} (120:1) -- (180:1) -- (0:1);
\draw[red,very thick] (120:1) -- (180:1);
    
    \node at (0,0.4) {$3$};
    
\end{scope}
\begin{scope}[xshift=1.5cm,yshift=-0.855cm]
    \draw (0:1) -- node[midway] {}(60:1) --  (120:1); 
    \draw[] (120:1) -- node[pos=0.5, left,xshift=2mm] {}(180:1);
   \draw (180:1) -- node[midway,below] {}(0:1);
\node at (-0.3,1.6) {$I$};    
    \node at (0,0.4) {$3$};
    
\end{scope}   
\end{scope}
\begin{scope}[xshift=5cm,yshift=-4.5cm]
     \draw (-1,0) .. controls (-6,1.5) and (0,2.3)  .. (2.5,0.85);
   \draw (0:1) -- node[midway,right] {}(60:1) -- (120:1) -- node[midway,left] {}(180:1);
    \draw (300:1) -- node[midway,below] {}(240:1);
\draw[red,very thick] (0:1) -- (60:1);
\draw[red,very thick] (240:1) -- (300:1);
\node at (0,0) {$2$};
\begin{scope}[xshift=-0.15cm, yshift=-0.5cm ,rotate=-30]
   \draw (-1,0) -- (-1,1) -- (-2,1) -- (-2,0) -- cycle;
   \draw[red,very thick] (-1,0) -- (-2,0);
   \draw[red,very thick] (-1,1) -- (-2,1);
     \node at (-1.5,0.5) {1}; 
\end{scope}
\begin{scope}[xshift=3cm]
\draw[] (0:1) -- node[pos=0.5, right, xshift=-4mm, yshift=0mm] {} (60:1);
    \draw (60:1) -- node[midway,above] {} (120:1) -- (180:1) --(240:1)-- (300:1) -- (0:1); 
    \node at (0,0) {$2$};
        \draw[red,very thick] (240:1) -- (180:1);
\draw[red,very thick] (0:1) -- (300:1);
 \draw[red,very thick] (120:1) -- (60:1);
\end{scope}

\begin{scope}[xshift=-1.5cm,yshift=-0.855cm]
    \draw[] (0:1) -- node[pos=0.5, right, xshift=-4mm, yshift=0mm] {} (60:1);
    \draw (60:1) -- node[midway,above] {} (120:1) -- (180:1) -- (0:1);
\draw[red,very thick] (120:1) -- (180:1);
    
    \node at (0,0.4) {$3$};
    
\end{scope}
\begin{scope}[xshift=1.5cm,yshift=-0.855cm]
    \draw (0:1) -- node[midway] {}(60:1) --  (120:1); 
    \draw[] (120:1) -- node[pos=0.5, left,xshift=2mm] {}(180:1);
   \draw (180:1) -- node[midway,below] {}(0:1);
\node at (-0.3,1.6) {$I$};    
    \node at (0,0.4) {$3$};
    
\end{scope}   
\end{scope}
\draw (-1,-7.5) --node[midway,left] {2} (-1,-6);
\draw (3,-7.3) --node[midway,above right] {2} (-1,-6);
\draw (3,-7.3) -- node[midway,right] {1}(3,-6);
\begin{scope}[xshift=-3cm,yshift=-9cm]
     \draw (-1,0) .. controls (-6,1.5) and (0,2.3)  .. (2.5,0.85);
   \draw (0:1) -- node[midway,right] {}(60:1) -- (120:1) -- node[midway,left] {}(180:1);
    \draw (300:1) -- node[midway,below] {}(240:1);
    \draw[red,very thick] (240:1) -- (180:1);
\draw[red,very thick] (0:1) -- (300:1);
 \draw[red,very thick] (120:1) -- (60:1);
\node at (0,0) {$2$};
\begin{scope}[xshift=-0.15cm, yshift=-0.5cm ,rotate=-30]
   \draw (-1,0) -- (-1,1) -- (-2,1) -- (-2,0) -- cycle;
   \draw[red,very thick] (-2,0) -- (-2,1);
     \node at (-1.5,0.5) {1}; 
\end{scope}
\begin{scope}[xshift=3cm]
\draw[] (0:1) -- node[pos=0.5, right, xshift=-4mm, yshift=0mm] {} (60:1);
    \draw (60:1) -- node[midway,above] {} (120:1) -- (180:1) --(240:1)-- (300:1) -- (0:1); 
    \node at (0,0) {$2$};
        \draw[red,very thick] (120:1) -- (180:1);
\draw[red,very thick] (240:1) -- (300:1);
 \draw[red,very thick] (0:1) -- (60:1);
\end{scope}

\begin{scope}[xshift=-1.5cm,yshift=-0.855cm]
    \draw[] (0:1) -- node[pos=0.5, right, xshift=-4mm, yshift=0mm] {} (60:1);
    \draw (60:1) -- node[midway,above] {} (120:1) -- (180:1) -- (0:1);
 \draw[red,very thick] (120:1) -- (180:1);
    
    \node at (0,0.4) {$3$};
    
\end{scope}
\begin{scope}[xshift=1.5cm,yshift=-0.855cm]
    \draw (0:1) -- node[midway] {}(60:1) --  (120:1); 
    \draw[] (120:1) -- node[pos=0.5, left,xshift=2mm] {}(180:1);
   \draw (180:1) -- node[midway,below] {}(0:1);
\node at (-0.3,1.6) {$I$};    
    \node at (0,0.4) {$3$};
    
\end{scope}   
\end{scope}
\begin{scope}[xshift=5cm,yshift=-9cm]
     \draw (-1,0) .. controls (-6,1.5) and (0,2.3)  .. (2.5,0.85);
   \draw (0:1) -- node[midway,right] {}(60:1) -- (120:1) -- node[midway,left] {}(180:1);
    \draw (300:1) -- node[midway,below] {}(240:1);
       \draw[red,very thick] (120:1) -- (180:1);
\draw[red,very thick] (240:1) -- (300:1);
 \draw[red,very thick] (0:1) -- (60:1);
\node at (0,0) {$2$};
\begin{scope}[xshift=-0.15cm, yshift=-0.5cm ,rotate=-30]
   \draw (-1,0) -- (-1,1) -- (-2,1) -- (-2,0) -- cycle;
   \draw[red,very thick] (-1,0) -- (-1,1);
   \draw[red,very thick] (-2,0) -- (-2,1);
     \node at (-1.5,0.5) {1}; 
\end{scope}
\begin{scope}[xshift=3cm]
\draw[] (0:1) -- node[pos=0.5, right, xshift=-4mm, yshift=0mm] {} (60:1);
    \draw (60:1) -- node[midway,above] {} (120:1) -- (180:1) --(240:1)-- (300:1) -- (0:1); 
    \node at (0,0) {$2$};
           \draw[red,very thick] (240:1) -- (180:1);
\draw[red,very thick] (0:1) -- (300:1);
 \draw[red,very thick] (120:1) -- (60:1);
\end{scope}

\begin{scope}[xshift=-1.5cm,yshift=-0.855cm]
    \draw[] (0:1) -- node[pos=0.5, right, xshift=-4mm, yshift=0mm] {} (60:1);
    \draw (60:1) -- node[midway,above] {} (120:1) -- (180:1) -- (0:1);

    \draw[red, very thick] (120:1) -- (180:1);
    \node at (0,0.4) {$3$};
    
\end{scope}
\begin{scope}[xshift=1.5cm,yshift=-0.855cm]
    \draw (0:1) -- node[midway] {}(60:1) --  (120:1); 
    \draw[] (120:1) -- node[pos=0.5, left,xshift=2mm] {}(180:1);
   \draw (180:1) -- node[midway,below] {}(0:1);
\node at (-0.3,1.6) {$I$};    
    \node at (0,0.4) {$3$};
    
\end{scope}   
\end{scope}
\begin{scope}[yshift=-4.5cm]
 \draw (-1,-7.5) --node[midway,left] {3} (-1,-6);
\draw (3,-7.3) -- node[midway,above right] {2}(-1,-6);
\draw (3,-7.3) -- node[midway,right] {2}(3,-6);   
\end{scope}

\begin{scope}[xshift=-3cm,yshift=-13.5cm]
     \draw (-1,0) .. controls (-6,1.5) and (0,2.3)  .. (2.5,0.85);
   \draw (0:1) -- node[midway,right] {}(60:1) -- (120:1) -- node[midway,left] {}(180:1);
    \draw (300:1) -- node[midway,below] {}(240:1);
\draw[red,very thick] (300:1) -- (0:1);
\draw[red,very thick] (120:1) -- (60:1);
\node at (0,0) {$2$};
\begin{scope}[xshift=-0.15cm, yshift=-0.5cm ,rotate=-30]
   \draw (-1,0) -- (-1,1) -- (-2,1) -- (-2,0) -- cycle;
   \draw[red,very thick] (-2,1) -- (-2,0);
     \node at (-1.5,0.5) {1}; 
\end{scope}
\begin{scope}[xshift=3cm]
\draw[] (0:1) -- node[pos=0.5, right, xshift=-4mm, yshift=0mm] {} (60:1);
    \draw (60:1) -- node[midway,above] {} (120:1) -- (180:1) --(240:1)-- (300:1) -- (0:1); 
    \node at (0,0) {$2$};
    \draw[red,very thick] (60:1) -- (0:1);
\draw[red,very thick] (120:1) -- (180:1);
\draw[red,very thick] (240:1) -- (300:1);
\end{scope}

\begin{scope}[xshift=-1.5cm,yshift=-0.855cm]
    \draw[] (0:1) -- node[pos=0.5, right, xshift=-4mm, yshift=0mm] {} (60:1);
    \draw (60:1) -- node[midway,above] {} (120:1) -- (180:1) -- (0:1);
\draw[red,very thick] (180:1) -- (0:1);
\draw[red,very thick] (120:1) -- (60:1);
    
    \node at (0,0.4) {$3$};
    
\end{scope}
\begin{scope}[xshift=1.5cm,yshift=-0.855cm]
    \draw (0:1) -- node[midway] {}(60:1) --  (120:1); 
    \draw[] (120:1) -- node[pos=0.5, left,xshift=2mm] {}(180:1);
   \draw (180:1) -- node[midway,below] {}(0:1);
\node at (-0.3,1.6) {$I$};    
    \node at (0,0.4) {$3$};
    
\end{scope}   
\end{scope}
\begin{scope}[xshift=5cm,yshift=-13.5cm]
     \draw (-1,0) .. controls (-6,1.5) and (0,2.3)  .. (2.5,0.85);
   \draw (0:1) -- node[midway,right] {}(60:1) -- (120:1) -- node[midway,left] {}(180:1);
    \draw (300:1) -- node[midway,below] {}(240:1);
    \draw[red,very thick] (60:1) -- (120:1);
\draw[red,very thick] (240:1) -- (180:1);
\draw[red,very thick] (0:1) -- (300:1);
\node at (0,0) {$2$};
\begin{scope}[xshift=-0.15cm, yshift=-0.5cm ,rotate=-30]
   \draw (-1,0) -- (-1,1) -- (-2,1) -- (-2,0) -- cycle;
   \draw[red,very thick] (-2,1) -- (-2,0);
     \node at (-1.5,0.5) {1}; 
\end{scope}
\begin{scope}[xshift=3cm]
\draw[] (0:1) -- node[pos=0.5, right, xshift=-4mm, yshift=0mm] {} (60:1);
    \draw (60:1) -- node[midway,above] {} (120:1) -- (180:1) --(240:1)-- (300:1) -- (0:1); 
    \node at (0,0) {$2$};
        \draw[red,very thick] (60:1) -- (120:1);
\draw[red,very thick] (240:1) -- (180:1);
\draw[red,very thick] (0:1) -- (300:1);
\end{scope}

\begin{scope}[xshift=-1.5cm,yshift=-0.855cm]
    \draw[] (0:1) -- node[pos=0.5, right, xshift=-4mm, yshift=0mm] {} (60:1);
    \draw (60:1) -- node[midway,above] {} (120:1) -- (180:1) -- (0:1);

        \draw[red,very thick] (180:1) -- (120:1);
    \node at (0,0.4) {$3$};
    
\end{scope}
\begin{scope}[xshift=1.5cm,yshift=-0.855cm]
    \draw (0:1) -- node[midway] {}(60:1) --  (120:1); 
    \draw[] (120:1) -- node[pos=0.5, left,xshift=2mm] {}(180:1);
   \draw (180:1) -- node[midway,below] {}(0:1);
\node at (-0.3,1.6) {$I$};    
    \node at (0,0.4) {$3$};
    
\end{scope}   
\end{scope}

\begin{scope}[yshift=-9cm]
 \draw (-1,-7.5) -- node[midway,left] {2} (-1,-6);
\draw (3,-7.3) -- node[midway,right] {3}(3,-6);   
\end{scope}

\begin{scope}[xshift=-3cm,yshift=-18cm]
     \draw (-1,0) .. controls (-6,1.5) and (0,2.3)  .. (2.5,0.85);
   \draw (0:1) -- node[midway,right] {}(60:1) -- (120:1) -- node[midway,left] {}(180:1);
    \draw (300:1) -- node[midway,below] {}(240:1);
    \draw[red,very thick] (60:1) -- (120:1);
\draw[red,very thick] (0:1) -- (300:1);
\node at (0,0) {$2$};
\begin{scope}[xshift=-0.15cm, yshift=-0.5cm ,rotate=-30]
   \draw (-1,0) -- (-1,1) -- (-2,1) -- (-2,0) -- cycle;
   \draw[red,very thick] (-2,1) -- (-2,0);
     \node at (-1.5,0.5) {1}; 
\end{scope}
\begin{scope}[xshift=3cm]
\draw[] (0:1) -- node[pos=0.5, right, xshift=-4mm, yshift=0mm] {} (60:1);
    \draw (60:1) -- node[midway,above] {} (120:1) -- (180:1) --(240:1)-- (300:1) -- (0:1); 
    \node at (0,0) {$2$};
        \draw[red,very thick] (60:1) -- (120:1);
\draw[red,very thick] (240:1) -- (180:1);
\draw[red,very thick] (0:1) -- (300:1);
\end{scope}

\begin{scope}[xshift=-1.5cm,yshift=-0.855cm]
    \draw[] (0:1) -- node[pos=0.5, right, xshift=-4mm, yshift=0mm] {} (60:1);
    \draw (60:1) -- node[midway,above] {} (120:1) -- (180:1) -- (0:1);
    \draw[red,very thick] (60:1) -- (120:1);
\draw[red,very thick] (0:1) -- (180:1);
    
    \node at (0,0.4) {$3$};
    
\end{scope}
\begin{scope}[xshift=1.5cm,yshift=-0.855cm]
    \draw (0:1) -- node[midway] {}(60:1) --  (120:1); 
    \draw[] (120:1) -- node[pos=0.5, left,xshift=2mm] {}(180:1);
   \draw (180:1) -- node[midway,below] {}(0:1);
\node at (-0.3,1.6) {$I$};    
    \node at (0,0.4) {$3$};
    
\end{scope}   
\end{scope}
\begin{scope}[xshift=5cm,yshift=-18cm]
     \draw (-1,0) .. controls (-6,1.5) and (0,2.3)  .. (2.5,0.85);
   \draw (0:1) -- node[midway,right] {}(60:1) -- (120:1) -- node[midway,left] {}(180:1);
    \draw (300:1) -- node[midway,below] {}(240:1);
    \draw[red,very thick] (60:1) -- (120:1);
\draw[red,very thick] (240:1) -- (180:1);
\node at (0,0) {$2$};
\begin{scope}[xshift=-0.15cm, yshift=-0.5cm ,rotate=-30]
   \draw (-1,0) -- (-1,1) -- (-2,1) -- (-2,0) -- cycle;
   \draw[red,very thick] (-2,1) -- (-2,0);
     \node at (-1.5,0.5) {1}; 
\end{scope}
\begin{scope}[xshift=3cm]
\draw[] (0:1) -- node[pos=0.5, right, xshift=-4mm, yshift=0mm] {} (60:1);
    \draw (60:1) -- node[midway,above] {} (120:1) -- (180:1) --(240:1)-- (300:1) -- (0:1); 
    \node at (0,0) {$2$};
        \draw[red,very thick] (60:1) -- (120:1);
\draw[red,very thick] (0:1) -- (300:1);
\end{scope}

\begin{scope}[xshift=-1.5cm,yshift=-0.855cm]
    \draw[] (0:1) -- node[pos=0.5, right, xshift=-4mm, yshift=0mm] {} (60:1);
    \draw (60:1) -- node[midway,above] {} (120:1) -- (180:1) -- (0:1);
    \draw[red,very thick] (180:1) -- (120:1);
    
    \node at (0,0.4) {$3$};
    
\end{scope}
\begin{scope}[xshift=1.5cm,yshift=-0.855cm]
    \draw (0:1) -- node[midway] {}(60:1) --  (120:1); 
    \draw[] (120:1) -- node[pos=0.5, left,xshift=2mm] {}(180:1);
   \draw (180:1) -- node[midway,below] {}(0:1);
\draw[red,very thick] (60:1) -- (120:1);
\draw[red,very thick] (180:1) -- (0:1);
\node at (-0.3,1.6) {$I$};
    \node at (0,0.4) {$3$};
    
\end{scope}   
\end{scope}
\begin{scope}[yshift=-23cm]
   \draw (-1,0) .. controls (-6,1.5) and (0,2.3)  .. (2.5,0.85);
   \draw (0:1) -- node[midway,right] {}(60:1) -- (120:1) -- node[midway,left] {}(180:1);
    \draw (300:1) -- node[midway,below] {}(240:1);
    \draw[red,very thick] (60:1) -- (120:1);
\node at (0,0) {$2$};
\begin{scope}[xshift=-0.15cm, yshift=-0.5cm ,rotate=-30]
   \draw (-1,0) -- (-1,1) -- (-2,1) -- (-2,0) -- cycle;
   \draw[red,very thick] (-2,1) -- (-2,0);
     \node at (-1.5,0.5) {1}; 
\end{scope}
\begin{scope}[xshift=3cm]
\draw[] (0:1) -- node[pos=0.5, right, xshift=-4mm, yshift=0mm] {} (60:1);
    \draw (60:1) -- node[midway,above] {} (120:1) -- (180:1) --(240:1)-- (300:1) -- (0:1); 
    \node at (0,0) {$2$};
        \draw[red,very thick] (60:1) -- (120:1);
\draw[red,very thick] (0:1) -- (300:1);
\end{scope}

\begin{scope}[xshift=-1.5cm,yshift=-0.855cm]
    \draw[] (0:1) -- node[pos=0.5, right, xshift=-4mm, yshift=0mm] {} (60:1);
    \draw (60:1) -- node[midway,above] {} (120:1) -- (180:1) -- (0:1);
 \draw[red,very thick] (60:1) -- (120:1);
\draw[red,very thick] (180:1) -- (0:1);
    
    \node at (0,0.4) {$3$};
    
\end{scope}
\begin{scope}[xshift=1.5cm,yshift=-0.855cm]
    \draw (0:1) -- node[midway] {}(60:1) --  (120:1); 
    \draw[] (120:1) -- node[pos=0.5, left,xshift=2mm] {}(180:1);
   \draw (180:1) -- node[midway,below] {}(0:1);
  \node at (-0.3,1.6) {$I$};  
    \node at (0,0.4) {$3$};
    \draw[red,very thick] (60:1) -- (120:1);
\draw[red,very thick] (180:1) -- (0:1);
\end{scope} 
\end{scope}

\begin{scope}[yshift=4cm]
   \draw[red, very thick] (-1,0) .. controls (-6,1.5) and (0,2.3)  .. (2.5,0.85);
   \draw (0:1) -- node[midway,right] {}(60:1) -- (120:1) -- node[midway,left] {}(180:1);
    \draw (300:1) -- node[midway,below] {}(240:1);
    \draw[red,very thick] (60:1) -- (120:1);
\draw[red,very thick] (240:1) -- (300:1);
\node at (0,0) {$2$};
\begin{scope}[xshift=-0.15cm, yshift=-0.5cm ,rotate=-30]
   \draw (-1,0) -- (-1,1) -- (-2,1) -- (-2,0) -- cycle;
   \draw[red,very thick] (-2,1) -- (-2,0);
     \node at (-1.5,0.5) {1}; 
\end{scope}
\begin{scope}[xshift=3cm]
\draw[] (0:1) -- node[pos=0.5, right, xshift=-4mm, yshift=0mm] {} (60:1);
    \draw (60:1) -- node[midway,above] {} (120:1) -- (180:1) --(240:1)-- (300:1) -- (0:1); 
    \node at (0,0) {$2$};
\draw[red,very thick] (240:1) -- (300:1);
\draw[red,very thick] (0:1) -- (60:1);
\end{scope}

\begin{scope}[xshift=-1.5cm,yshift=-0.855cm]
    \draw[] (0:1) -- node[pos=0.5, right, xshift=-4mm, yshift=0mm] {} (60:1);
    \draw (60:1) -- node[midway,above] {} (120:1) -- (180:1) -- (0:1);
\draw[red,very thick] (180:1) -- (120:1);
    
    \node at (0,0.4) {$3$};
    
\end{scope}
\begin{scope}[xshift=1.5cm,yshift=-0.855cm]
    \draw (0:1) -- node[midway] {}(60:1) --  (120:1); 
    \draw[] (120:1) -- node[pos=0.5, left,xshift=2mm] {}(180:1);
   \draw (180:1) -- node[midway,below] {}(0:1);
    
    \node at (0,0.4) {$3$};
    \node at (-0.3,1.6) {$I$};
    \draw[red,very thick] (60:1) -- (120:1);
\end{scope} 
\end{scope}

\begin{scope}[yshift=-20.5cm]
 \draw (-1,1.5) -- node[midway,left] {3}(1,-1);
\draw (3,1.5) -- node[midway,right] {3}(1,-1);  
\end{scope}
   \end{tikzpicture}
   \caption{The poset of perfect matchings of $\mathcal{G}_{ab}$, and the corresponding monomials.}
   \label{set_pm}
\end{figure} 

\section*{Acknowledgments}
The results presented in this article are part of my Ph.D. thesis under the supervision of Giovanni Cerulli Irelli. I thank him and Ralf Schiffler for many helpful discussions on the subject. I also thank Javier De Loera for introducing me to the lattice structure for orientations of graphs. Finally, I am grateful to the anonymous referees for their valuable comments.

This work has been partially supported by the "National Group for Algebraic and Geometric Structures, and their Applications" (GNSAGA - INdAM), by the National PRIN project "SQUARE", by Sapienza Progetti Grandi di Ateneo 2023 "Representation Theory and Applications", by Sapienza Progetti Piccoli 2022 "Algebraic varieties associated to finite dimensional algebras" and by Sapienza Progetti di Avvio alla Ricerca 2022 "Homological interpretation of cluster algebras of finite type through symplectic and orthogonal representations of symmetric
quivers".

\printbibliography[heading=bibintoc]

\end{document}